\pdfoutput=1
\documentclass[11pt,a4paper]{article}
\usepackage[utf8]{inputenc}
\usepackage[T1]{fontenc}
\usepackage{lmodern}
\usepackage[margin=2.5cm]{geometry}
\usepackage{amsmath,amssymb,amsthm}
\usepackage{graphicx}
\usepackage{booktabs}
\usepackage{array}
\usepackage{algorithm}
\usepackage{algpseudocode}
\usepackage[numbers,sort&compress]{natbib}
\usepackage{xcolor}
\usepackage{xurl}
\input{glyphtounicode}\pdfgentounicode=1
\usepackage[colorlinks=true,linkcolor=blue!50!black,citecolor=blue!50!black,urlcolor=blue!50!black]{hyperref}
\usepackage{caption}
\usepackage{subcaption}
\usepackage{adjustbox}
\usepackage{placeins}
\usepackage{enumitem}
\graphicspath{{figures/}}
\hypersetup{pdftitle={Territorial Conflict as a Cellular Automaton: Frozen Domains and Active Fronts in a Risk-Derived Local Rule},pdfauthor={María Teresa Arias Noguerales}}
\newtheorem{definition}{Definition}
\newtheorem{proposition}{Proposition}
\newtheorem{lemma}{Lemma}
\newtheorem{hypothesis}{Hypothesis}

\newcommand{\Nmax}{N_{\max}}

\begin{document}
\begin{center}
{\Large\bfseries Territorial Conflict as a Cellular Automaton:\\[2pt] Frozen Domains and Active Fronts in a \emph{Risk}-Derived Local Rule\par}
\vspace{1em}
{\large María Teresa Arias Noguerales\par}
\end{center}

\begin{abstract}
\noindent
We formulate the attack and fortification rules of the game \emph{Risk} as a synchronous cellular automaton on a graph, with battles resolved by exact dice odds. We prove, and verify in Lean, that for the default exchange ratio the rule has radius 6, or 5 for deterministic battles with a positive attack threshold, both attained; that its deterministic form is irreversible except at the highest attack thresholds without growth and with fortification absent or inactive; that without troop growth it freezes in finite time; and that the dynamics near an interface is locally determined. The coupling of the ownership field to a troop field separates frozen domains from active fronts: with a positive attack threshold, saturated domains cannot be conquered, while interfaces with heterogeneous troops keep fighting. In 14\,252 runs, with attractors classified operationally in Wolfram's scheme, raising the attack threshold takes the system from dominance by one faction to a frozen fragmented partition. The share of the largest faction drops between the same two threshold values at every size up to $L=128$, in a step about 0.5 wide under a continuous threshold ($L=64$--96), whereas the probability of complete consolidation is a size-dependent crossover governed by saturated minority domains. Long-period attractors combine spatially fixed local oscillators; apparent class-IV gliders are artefacts of index-ordered tie-breaking. Straight fronts stayed within two cells of their initial position for up to $5\times10^5$ steps, over which they are one-dimensional automata with non-recurrent, incompressible dynamics and positive Lyapunov speeds, constant for $L\ge48$. Compared with mean-field and random-regular-graph controls, the lattice lengthens conflicts about twentyfold and clusters conquest: a cell next to a fresh conquest is about five times as likely as average to be conquered next, against 1.1--1.5 in the controls.
\end{abstract}


\section{Introduction}
\label{sec:intro}

\emph{Risk} models war with territories, an adjacency relation, integer armies, random attrition, reinforcement and redeployment, the ingredients of agent-based models of combat and civil conflict \cite{ilachinski2004,epstein2002}. Since the board is a graph whose nodes carry a finite state updated through their neighbours, the combat and fortification rules of the game can be written as a cellular automaton (CA) on a graph \cite{wolfram1983}.

Mathematical work on \emph{Risk} has not used this description. It treats a single battle as an absorbing Markov chain \cite{tan1997,osborne2003}, which gives exact odds but no spatial information. Lattice and agent-based models of combat \cite{ilachinski2004}, state formation \cite{cederman1997}, spatial competition \cite{clifford1973} and insurgency \cite{pratt2012} have, to our knowledge, not been applied to a \emph{Risk}-derived rule, and none places the dynamics of such a rule in one of Wolfram's classes \cite{wolfram1984}.

A lattice formulation describes what the battle chain cannot: the spatial spread of conquest, the formation and coarsening of fronts, and the structure of attractors. The model studied here, RiskCA, is derived from \emph{Risk} but is not the game. Every cell follows one fixed heuristic, ``attack the weakest hostile neighbour when stronger by $\theta$''; a fixed growth rate replaces reinforcement; each cell holds at most $\Nmax$ units; only the battle kernel is taken from the game unchanged. RiskCA is also far from the elementary and totalistic rules for which Wolfram's classes were introduced: each cell has $P\Nmax$ states (64 by default) and the rule has radius up to 6 (Lemma~\ref{lem:radius}). We study one family of rules, parametrised by growth and attack margin, and use the classification to describe its dynamics, not to test the scheme or to predict the course of games or conflicts. We ask:
\begin{enumerate}[label=(Q\arabic*),leftmargin=*,itemsep=0pt]
\item What formal properties does the rule have?
\item How does the attack margin $\theta$ determine the long-run behaviour, and is the change from dominance to fragmentation a phase transition?
\item Which observed behaviours are long-period or long-transient class-II dynamics, which are sustained irregular dynamics, and which are artefacts of the implementation?
\item How does spatial structure change the dynamics, compared with a mean-field control and a random graph of the same degree?
\end{enumerate}

In brief, the answers are as follows.
\begin{itemize}[leftmargin=*,itemsep=0pt]
\item[(Q1)] The rule has radius up to 6, a bound that is attained, its deterministic form is irreversible for the default parameters, and without growth it always freezes; near an interface its dynamics is locally determined (Section~\ref{sec:model}). All these results are machine-checked in Lean~4 (Section~\ref{sec:design}).
\item[(Q2)] Low $\theta$ leads to dominance by one faction, high $\theta$ to a frozen partition. For $g=1$, the largest faction's share drops between $\theta=4$ and 5, more steeply in larger systems; with a continuous threshold the drop is a step about half a unit of $\theta$ wide at $\theta\approx4.7$, but whether it is a phase transition cannot be decided with our data (Section~\ref{sec:res:transition}).
\item[(Q3)] Resolved runs end on homogeneous or heterogeneous attractors, labelled I and II in the operational sense of Section~\ref{sec:criteria}. The informative findings are about their structure: long cycles combine spatially fixed local oscillators, apparent class-IV structures are tie-breaking artefacts, and sustained irregular dynamics was observed only on stalemated fronts (Sections~\ref{sec:res:artefact} and \ref{sec:res:robust}).
\item[(Q4)] On the lattice, conflicts last about twenty times longer and conquests are more clustered than in either control (Section~\ref{sec:res:spatial}).
\end{itemize}

The main contribution is a mechanism rather than a classification. Several features of RiskCA, such as $L^2$ coarsening and stripe states, are shared with majority-type lattice models; what is specific is the coupling of the ownership field to a troop field, which separates the global attractor dynamics from irregular dynamics localised at interfaces: saturated troops make minority domains unconquerable and freeze fragmented partitions, while heterogeneous troops keep interfaces active where the corresponding Ising stripes are frozen. In our runs irregular dynamics is confined to one-dimensional fronts, and the comparison with a random regular graph points to the geometry of the lattice, rather than the locality of interaction, as what makes conflicts long and conquest clustered.

Section~\ref{sec:model} defines RiskCA and answers (Q1). Section~\ref{sec:results} addresses (Q2)--(Q4) and simulates the 42-territory board. Sections~\ref{sec:discussion} and \ref{sec:limits} discuss the results and their limitations.

\section{Related work}
\label{sec:related}

\paragraph{Complexity classes.} \citet{wolfram1983,wolfram1984} divided the long-run behaviour of CA into homogeneous states (class I), separated simple or periodic structures (II), aperiodic patterns (III), and interacting localised structures with long transients (IV). \citet{packard1985} extended the classification to two-dimensional rules, mainly two-state totalistic and outer-totalistic rules, and found the same four classes on the whole lattice; \citet{wolfram2002} surveys the classes across many kinds of systems. Class membership is undecidable in general \cite{culik1988}, so classifications rest on operational criteria such as Langton's $\lambda$ parameter \cite{langton1990}, damage spreading \cite{derrida1986,grassberger1995}, input-entropy variance \cite{wuensche1999} and compressibility \cite{zenil2010}; we use the last three. Invertible CA are reviewed in \cite{toffoli1990} and the effect of the updating scheme in \cite{schonfisch1999}; on irregular graphs the topology itself changes the patterns a rule can produce \cite{marr2005}. For one-dimensional CA, the speeds at which perturbations spread to the left and right define Lyapunov exponents \cite{shereshevsky1992}.

\paragraph{Spatial conflict and coarsening.} The voter model, introduced as ``a model for spatial conflict'' \cite{clifford1973}, has no surface tension: on the two-dimensional torus its consensus time grows as $N\log N$ \cite{cox1989} and its interface density decays logarithmically \cite{dornic2001}. With surface tension, as in majority-vote models \cite{deoliveira1992} and zero-temperature Ising dynamics, domains coarsen with a characteristic length growing as $t^{1/2}$ \cite{bray1994}. On the square lattice, zero-temperature Ising dynamics does not always reach consensus: a finite fraction of runs end in frozen states of straight stripes spanning the system \cite{spirin2001}. In spatial games, local imitation of more successful neighbours produces persistent and chaotically fluctuating spatial patterns \cite{nowak1992}.

RiskCA shares curvature-driven coarsening and stripe states with these models but not the mechanism: domains change through battles decided by integer troop counts, which growth raises, combat depletes and fortification redistributes (Section~\ref{sec:discussion}).

\paragraph{Analyses of Risk and models of combat.} \citet{tan1997} treated a battle as an absorbing Markov chain, and \citet{osborne2003} corrected the independence assumption for sorted dice. ISAAC and EINSTein \cite{ilachinski2004} analyse emergent behaviour in agent-based combat, \citet{epstein2002} models civil violence, \citet{pratt2012} couple fuzzy cognitive maps to a CA to model insurgency, and GeoSim \cite{cederman1997} models state formation through conquest on a grid. To our knowledge, none of these studies the collective dynamics of a \emph{Risk}-derived rule or places such dynamics in Wolfram's classes (Table~\ref{tab:positioning}).

\begin{table}[t]
\centering\small
\caption{Coverage of related work. \checkmark: addressed; ($\checkmark$): partly, e.g.\ agents moving on a lattice rather than a cellular automaton, or complex dynamics discussed without assigning Wolfram classes; --: not addressed.}
\label{tab:positioning}
\begin{adjustbox}{max width=\textwidth}
\begin{tabular}{lcccccc}
\toprule
 & Spatial & Exact battle & Local rule & Wolfram & Damage / & Board-level \\
 & lattice & kernel & as CA & class & attractors & statistics\\
\midrule
Markov analyses of Risk \cite{tan1997,osborne2003} & -- & \checkmark & -- & -- & -- & -- \\
ISAAC / EINSTein \cite{ilachinski2004} & \checkmark & -- & ($\checkmark$) & ($\checkmark$) & ($\checkmark$) & -- \\
FCM + CA insurgency \cite{pratt2012} & \checkmark & -- & \checkmark & -- & -- & -- \\
GeoSim \cite{cederman1997} & \checkmark & -- & ($\checkmark$) & -- & -- & -- \\
Voter, majority and Ising models \cite{clifford1973,dornic2001,deoliveira1992,spirin2001} & \checkmark & -- & \checkmark & ($\checkmark$) & \checkmark & -- \\
Spatial games \cite{nowak1992} & \checkmark & -- & \checkmark & ($\checkmark$) & -- & -- \\
This paper & \checkmark & \checkmark & \checkmark & \checkmark & \checkmark & \checkmark \\
\bottomrule
\end{tabular}
\end{adjustbox}
\end{table}

\section{The RiskCA model}
\label{sec:model}
The main symbols are listed in Appendix~\ref{app:notation}.

\subsection{Lattice and state}

\begin{definition}[Graph lattice]
A \emph{board} is a finite simple graph $\mathcal{G}=(V,E)$ with $|V|=N$ cells. The neighbourhood of cell $i$ is $\mathcal{N}(i)=\{j:(i,j)\in E\}$ and its degree is $k_i=|\mathcal{N}(i)|$.
\end{definition}

We use the $L\times L$ torus $\mathbb{T}_L$ with the von Neumann neighbourhood ($k_i=4$), and in Section~\ref{sec:res:robust} also open boundaries and the Moore neighbourhood. The other board is the 42-territory, six-region map of the game (83 edges, degrees 2--6, region bonuses $\{5,2,5,3,7,2\}$; Appendix~\ref{app:board}).

\begin{definition}[Cell state]
The state of cell $i$ at time $t$ is $x_i(t)=(o_i(t),n_i(t))\in\Sigma=\{1,\dots,P\}\times\{1,\dots,\Nmax\}$, where $o_i$ is the owning faction and $n_i$ the number of units stationed there. A configuration is $x\in\Sigma^V$; the configuration space has $(P\Nmax)^N$ elements.
\end{definition}

Every cell is always owned and holds at least one unit; the cap $\Nmax$ keeps the state space finite.

\subsection{Local update rule}
\begin{figure}[t]
\centering
\includegraphics[width=\textwidth]{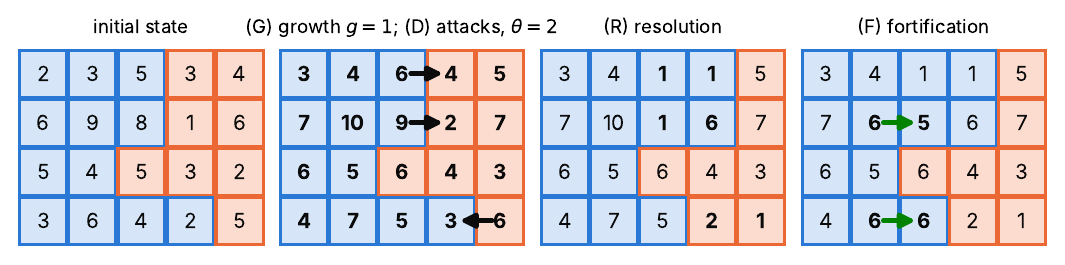}
\caption{One RiskCA step on a $4\times5$ open patch (two factions, $g=1$, $\theta=2$, deterministic battles), computed with the model code. Bold numbers mark cells changed by each sub-rule. (D) Arrows point from attackers to their targets. (R) Attackers keep one unit; all three attacks succeed. (F) Two interior cells send half their surplus towards the front.}
\label{fig:schematic}
\end{figure}

A time step applies $\Phi=\mathcal{F}\circ\mathcal{R}\circ\mathcal{D}\circ\mathcal{G}$, four sub-rules acting simultaneously on all cells (Figure~\ref{fig:schematic}). Write $H_i=\{j\in\mathcal{N}(i): o_j\neq o_i\}$ for the \emph{hostile} neighbours of $i$ and $\eta_i=|H_i|$ for its \emph{threat}. A cell with $\eta_i=0$ is \emph{interior}, otherwise \emph{frontier}. Ties are broken by a fixed priority $\pi$, a random permutation of $V$ drawn once per run (``quenched''); Section~\ref{sec:res:artefact} compares this with index order, $\pi(i)=i$. Because $\pi$ is attached to cells, a run with fixed $\pi$ is a CA with quenched spatial disorder in its tie-breaking: the update of a cell depends on the priorities around it, and the rule is not translation invariant. Index order is not translation invariant either; it ranks cells along a fixed direction, which is why it produces translating structures (Section~\ref{sec:res:artefact}). Each run draws its own $\pi$, so results with quenched priorities refer to the ensemble of disordered realisations.

\paragraph{(G) Growth.} $n_i\leftarrow\min(\Nmax,\,n_i+g)$ with integer $g\ge 0$. On the 42-territory board we instead use \emph{Risk-like} reinforcement: faction $p$ receives $R_p=\max(3,\lfloor|V_p|/3\rfloor)+\sum_{c\in\mathcal{C}_p}b_c$ units ($V_p$: cells owned by $p$; $\mathcal{C}_p$: regions it owns completely; $b_c$: region bonus), placed at random on its frontier cells (with the cap $\Nmax$). Because $R_p$ depends on $|V_p|$, this couples distant cells: the board model is a territorial dynamics with global feedback, not a CA. All formal results below, and all experiments except E7, concern the local rule with growth $g$.

\paragraph{(D) Decision.} Each frontier cell selects the weakest hostile neighbour, with ties broken by $\pi$ (pairs are compared lexicographically, here and below),
\begin{equation}
  j^*(i)=\arg\min_{j\in H_i}\,(n_j,\pi(j)),
\end{equation}
and \emph{wants} to attack it iff
\begin{equation}
  n_i\ge 2\quad\text{and}\quad n_i-n_{j^*(i)}\ge\theta .
  \label{eq:attack}
\end{equation}
The integer \emph{attack margin} $\theta\ge 0$ is the only parameter of the heuristic (E17 interpolates between integer values); since $n_i-n_j\le\Nmax-1$, no attacks occur for $\theta\ge\Nmax$.

\paragraph{(R) Resolution.} Each target $j$ accepts, among the cells that want to attack it, the one with the largest committed force $a_i=n_i-1$ (ties by $\pi$); the others hold, and a target that nobody wants to attack is not attacked. All accepted battles are resolved simultaneously from the configuration after (D). The attacker keeps one unit ($n_i\leftarrow 1$); the defender fights with its strength after commitments, $d=n_j$, so a cell whose own attack was accepted defends with one unit. A conquered cell takes the owner that its attacker had at the start of the step, even if the attacker's own cell is conquered in the same step; for example, if $A$ attacks $B$ and $B$ attacks $C$, both successfully, then $C$ passes to $B$'s original faction and $B$ to $A$'s. Two battle operators are used.
\begin{itemize}[leftmargin=*,itemsep=0pt]
\item \emph{Deterministic} ($\mathcal{R}_\rho$). The attacker wins iff $\rho a>d$ and moves in with $a'=\max(1,a-\lceil d/\rho\rceil)$ units; otherwise the defender keeps $d'=\max(1,d-\lfloor\rho a\rfloor)$. The exchange ratio $\rho=1.1718$ (1.171764 in the code) is the ratio of expected defender to attacker losses in a 3-against-2 roll (Section~\ref{sec:res:dice}), so $\mathcal{R}_\rho$ replaces the dice by their expected losses. It reproduces the mean exchange of units but not the variance of battle outcomes; how far this matters is tested against the exact dice in Sections~\ref{sec:res:dice} and \ref{sec:discussion}.
\item \emph{Exact dice} ($\mathcal{R}_{\mathrm{dice}}$). The outcome is drawn from the absorbing distribution $K(a,d;\cdot)$ of the battle chain \cite{tan1997,osborne2003}: the attacker rolls $\min(3,a)$ dice and the defender $\min(2,d)$, dice are sorted and paired, ties go to the defender, and the battle continues until one side is exhausted. $K$ is computed by dynamic programming (Algorithm~\ref{alg:kernel}).
\end{itemize}

\paragraph{(F) Fortification.} After battle, each interior cell $i$ sends $m_i=\lfloor(n_i-n_{k^*})/2\rfloor$ units to the friendly neighbour of highest threat, $k^*(i)=\arg\max_{k\in\mathcal{N}(i),\,o_k=o_i}(\eta_k,-\pi(k))$, provided $\eta_{k^*}>0$ and $m_i\ge1$. Inflows are capped at $\Nmax$ and the excess (\emph{overflow}) is lost. Receivers have $\eta>0$, so no cell both sends and receives.

\paragraph{Schedules and control.} In the \emph{synchronous} schedule, which defines the CA, all cells apply $\Phi$ at once. In the \emph{turn} schedule, closer to the game, each round has $P$ phases in random order, and only the cells of faction $p$ act in phase $p$.

In the \emph{mean-field} control, $\mathcal{N}(i)$ in (D)--(R) is replaced by $k_i$ cells other than $i$, drawn uniformly with replacement for each cell and redrawn at every step (for $N\ge256$, fewer than 3\% of cells draw the same partner twice); growth and fortification keep the lattice neighbourhood. This removes spatial correlations in combat but also changes who meets whom: a cell faces the global rather than the local distribution of owners and troops. A second control runs the unchanged rule on a quenched random 4-regular graph, which keeps the degree and the locality of interactions but has no low-dimensional geometry.

\begin{lemma}[Radius]
\label{lem:radius}
Let the radius of the rule be the smallest $r$ such that the state of every cell after one step depends only on the states of the cells within graph distance $r$ of it. For dice, the step is a function of the configuration and the random numbers; the radius refers to this function with the random numbers held fixed, so a bound holds for every realisation and is attained if it is attained for some realisation. On any board and for any priority $\pi$:
(a) the radius is at most 6;
(b) with deterministic battles, $\rho>1$ and $\theta\ge1$, every accepted attack succeeds and the radius is at most 5;
(c) for $\theta\ge\Nmax$, when no attack is possible, the radius is at most 3.
On open patches of the square lattice the bound 6 is attained with deterministic battles and $\theta=0$, and with dice and $\theta=1$; the bound 5 is attained with deterministic battles and $\theta=1$.
\end{lemma}
\begin{proof}
(a) A decision in (D) depends on the cell and its neighbours (distance 1). Acceptance at a target $j$ depends on the decisions of the neighbours of $j$, hence on states within distance 2 of $j$. The state of a cell after (R) depends on whether its own attack was accepted (distance $1+2$), on whether and by whom it was conquered (distance 2), and on the committed forces involved; it therefore depends on states within distance 3. In (F), a sender $i$ needs the owners within distance 2 (to compute $\eta_i$ and the threats of its neighbours) and the troops of its destination; a receiver adds the inflows of its neighbours. The result of (F) at a cell therefore depends on the post-(R) states within distance 3, so $\Phi$ has radius at most $3+3=6$.
(b) An attack on $j$ is wanted only if $n_i-n_j\ge\theta\ge1$, so $a=n_i-1\ge n_j\ge d$, and $\rho>1$ gives $\rho a>d$. Whether a target is conquered, and by whom, then depends only on acceptance at that target, hence on states within distance 2 of it. Since owners change only at targets, the post-(R) owner of any cell depends on states within distance 2; in particular, whether an attacker's own cell is conquered depends on acceptance at that cell and not on whether its own attack elsewhere was accepted. The post-(R) troops of a cell depend in addition on whether its own attack was accepted, hence on states within distance 3. A sender's choice and transfer then depend on states within distance 4, and a receiver's inflow on states within distance 5.
(c) Without attacks the owners do not change, a sender depends on the owners within distance 2 of itself and on the troops of itself and its destination, and a receiver on its neighbouring senders, hence on states within distance 3.
The bound in (a) is independent of $\pi$, which enters only through tie-breaking. Appendix~\ref{app:witness} gives an explicit pair of configurations attaining 6 with deterministic battles and $\theta=0$. Configurations attaining 6 with dice and $\theta=1$ (for fixed random numbers) and 5 with deterministic battles and $\theta=1$ were found by local search on a $7\times6$ patch; all three are reproduced by unit tests.
\end{proof}
A perturbation can spread at most 6 cells per step, and at most 5 with deterministic battles, $\rho>1$ and $\theta\ge1$.

\subsection{Formal properties}

\begin{proposition}[Finite-state determinism]
\label{prop:finite}
For the synchronous deterministic rule ($\mathcal{R}_\rho$, integer $g$, fixed $\pi$), every orbit is eventually periodic: there exist a transient $\tau\ge 0$ and a period $p\ge1$ with $x(\tau+p)=x(\tau)$.
\end{proposition}
\begin{proof}
$\Phi$ is a function on the finite set $\Sigma^V$. By the pigeonhole principle, some configuration repeats within $(P\Nmax)^N+1$ steps.
\end{proof}
With $\mathcal{R}_{\mathrm{dice}}$ the dynamics is a Markov chain on $\Sigma^V$, and individual trajectories need not be periodic. Throughout the paper, transients, periods, attractors and Wolfram classes refer only to the deterministic rule. For the dice rule we use only absorbing configurations, in which no attack is possible so that the step is deterministic (for example consolidated or saturated configurations), and observed quantities such as the time to consolidation and the rate of owner changes. Propositions~\ref{prop:irrev} and \ref{prop:fixed}(c) also assume deterministic battles; Propositions~\ref{prop:balance}, \ref{prop:fixed}(a,b), \ref{prop:g0} and \ref{prop:local} hold for both operators.

\begin{proposition}[Irreversibility]
\label{prop:irrev}
Let the board be connected with $N\ge2$, let $\Nmax\ge3$, and let battles be deterministic. Then $\Phi$ is not injective if at least one of the following holds: (i) $g\ge1$; (ii) $\theta\le1$ and $\rho>1$; (iii) $\theta\ge2$, fortification is enabled and the board is not a complete graph; (iv) $2\le\theta\le\Nmax-2$ and $1<\rho<2$. In each case there are configurations without a pre-image (Garden-of-Eden configurations). Without growth, and without fortification or on a complete graph, where (F) is inactive, two cases remain: if $\theta\ge\Nmax$, no attack is possible and $\Phi$ is the identity, hence injective; if $\theta=\Nmax-1$, attacks are possible only from saturated cells against cells with one unit, and injectivity is left open. For $1<\rho<2$, which includes the default, these are the only cases not covered.
\end{proposition}
\begin{proof}
The conditions can hold simultaneously; since (i) covers $g\ge1$, we assume $g=0$ in (ii)--(iv). In each case we give two distinct configurations $X\neq Y$ with $\Phi(X)=\Phi(Y)$.
(i) Let $X$ be homogeneous with $n_i\equiv\Nmax$, and $Y$ equal to $X$ except $n_u=\Nmax-1$ at one cell. No cell has a hostile neighbour, so (D), (R) and (F) are inactive, and (G) maps both to $X$.
(ii) Take an edge $(u,v)$. In $X$, cell $v$ is owned by faction $B$ with $n_v=1$, cell $u$ by $A$ with $n_u=2$, and every other cell by $A$ with one unit; $Y$ equals $X$ except $n_u=3$. Only $u$ can attack, since all other cells hold one unit, and it attacks $v$ because $n_u-n_v\ge1\ge\theta$. With $a=1$ or $a=2$ against $d=1$ the attacker wins ($\rho>1$) and keeps $\max(1,a-\lceil1/\rho\rceil)=1$ unit. Both images are homogeneous with one unit per cell, and (F) is inactive.
(iii) Since the board is connected and not complete, it has cells $u,v$ at graph distance 2. Give $v$ to faction $B$ and all other cells to $A$, so that $u$ is interior and has a friendly neighbour adjacent to $v$; let $k=k^*(u)$, which depends only on ownership and has $\eta_k\ge1$. In $X$ set $n_u=3$, $n_k=1$, in $Y$ set $n_u=2$, $n_k=2$, and give every other cell one unit. No battle occurs, since every cell except $u$ and $k$ holds one unit, $u$ is interior and $n_k-n_v\le1<\theta$. In (F), only $u$ can send; it transfers $\lfloor(3-1)/2\rfloor=1$ unit in $X$ and none in $Y$, so both images have $n_u=n_k=2$. This case does not involve $\rho$.
(iv) Take an edge $(u,v)$, give $v$ to faction $B$ and every other cell to $A$ with one unit. In $X$ set $n_u=\theta+1$, $n_v=1$; in $Y$ set $n_u=\theta+2$, $n_v=2$, which is possible since $\theta+2\le\Nmax$. Only $u$ can attack: other cells of $A$ hold one unit, and $n_v-1\le1<\theta$. Since $n_u-n_v=\theta$, $u$ attacks $v$ with $a=\theta$ against $d=1$ in $X$ and $a=\theta+1$ against $d=2$ in $Y$, and wins in both because $a>d$ and $\rho>1$. For $1<\rho<2$, $\lceil1/\rho\rceil=1$ and $\lceil2/\rho\rceil=2$, so the conquering force is $\theta-1$ in both cases. Both images are owned by $A$ throughout, with $n_u=1$, $n_v=\theta-1$ and one unit elsewhere, and (F) is inactive.
Since a map from a finite set to itself is surjective iff it is injective, $\Phi$ has Garden-of-Eden configurations.
\end{proof}

\begin{proposition}[Troop balance]
\label{prop:balance}
Let $M(t)=\sum_i n_i(t)$ be the troop mass and $B(t)$ the number of battles at step $t$. Then
\begin{equation}
M(t+1)-M(t)=\Gamma(t)-\Lambda(t)-\Omega(t),\qquad 0\le\Gamma(t)\le gN,\quad \Omega(t)\ge0,\quad \Lambda(t)\ge B(t),
\end{equation}
where $\Gamma$ is the growth actually added (after capping), $\Lambda$ the units destroyed in combat and $\Omega$ the units lost to the cap in fortification. Moreover $N\le M(t)\le N\Nmax$.
\end{proposition}
\begin{proof}
(G) adds $\Gamma$ units. (F) conserves units except for the capped excess $\Omega$. For a battle with committed force $a\ge1$ and defence $d\ge1$, the attacker side loses $a$ units from its origin and the pair ends with $a'$ (win) or $d'$ (loss). Under $\mathcal{R}_\rho$ a win loses $a+d-a'\ge d\ge1$ units because $a'\le a$, and a defeat loses $a+d-d'\ge a\ge1$ because $d'\le d$. Under $\mathcal{R}_{\mathrm{dice}}$ the same inequalities hold for the survivor counts. A cell whose own attack is accepted and which is attacked in the same step takes part in two battles, but no unit is counted twice: its committed force $a$ belongs to the first battle and the single unit it keeps is its defence $d=1$ in the second, so $\Lambda(t)$ is the sum of the losses of the $B(t)$ battles, each at least one. The bounds follow from $1\le n_i\le\Nmax$.
\end{proof}
Troop mass is not conserved: growth injects it, and combat and overflow dissipate it.

\begin{proposition}[Fixed points and absorbing set]
\label{prop:fixed}
(a) The homogeneous configurations ($o_i\equiv p$) form an absorbing set. For $g\ge1$ they reach the saturated homogeneous fixed point within $\lceil(\Nmax-1)/g\rceil$ steps.
(b) Let $\theta\ge1$. Then, for any $g$, every \emph{saturated} configuration ($n_i\equiv\Nmax$), with any ownership pattern, is a fixed point; there are $P^N$ of them.
(c) Let the board be connected, battles deterministic, $g\ge1$, $\theta=0$ and $\rho(\Nmax-1)>\Nmax$ (true for the defaults $\rho=1.1718$, $\Nmax=16$). Then a saturated configuration is a fixed point iff it is homogeneous.
\end{proposition}
\begin{proof}
(a) Without hostile neighbours, (D)--(F) are inactive and (G) raises every cell by $g$ until the cap. (b) Saturation makes (G) a no-op. Condition~\eqref{eq:attack} requires $n_i-n_j\ge\theta\ge1$, but all differences are zero, and (F) transfers $\lfloor 0/2\rfloor=0$. (c) On a connected board a heterogeneous configuration has at least one frontier cell, and when it is saturated every frontier cell satisfies \eqref{eq:attack} with $\theta=0$. Each selected attack commits $a=\Nmax-1$ against $d\le\Nmax$ and wins, because $\rho(\Nmax-1)>\Nmax\ge d$, so ownership changes.
\end{proof}
If $\rho(\Nmax-1)\le\Nmax$, an attack on a saturated defender that has not committed itself fails, ownership may be unchanged, and whether a heterogeneous saturated configuration is a fixed point depends on whether fortification restores the troops lost in combat; we have not characterised this case, which does not arise for the default parameters.

\begin{proposition}[Finite-time freezing without growth]
\label{prop:g0}
If $g=0$, every trajectory, with deterministic or dice battles, reaches a fixed point after finitely many steps; the number of steps in which the configuration changes is at most $B_{\max}(1+N\Nmax)+N\Nmax$, with $B_{\max}=N(\Nmax-1)$.
\end{proposition}
\begin{proof}
With $g=0$, (G) is inactive. Let $K=N\Nmax$ and $V=(K+1)M+M_{\mathrm{int}}$, where $M_{\mathrm{int}}=\sum_{i:\eta_i=0}n_i$ is the troop mass on interior cells. Since $0\le M_{\mathrm{int}}\le K$ and $N\le M\le K$, we have $(K+1)N\le V\le(K+1)K+K$. We show that $V$ decreases by at least one in every step that changes the configuration. In a step with a battle, Proposition~\ref{prop:balance} gives $M(t+1)\le M(t)-1$, which lowers $(K+1)M$ by at least $K+1$, more than any increase of $M_{\mathrm{int}}$. In a step without battles, (R) changes nothing, so ownership and the threats $\eta_i$ are unchanged and only (F) acts. It moves units from interior cells to frontier cells and never in the opposite direction, because receivers have $\eta>0$; hence $M$ does not increase and $M_{\mathrm{int}}$ decreases by the number of units sent, which is at least one if the configuration changes. The number of changing steps is therefore at most $(K+1)K+K-(K+1)N=B_{\max}(1+K)+K$. If a step does not change the configuration, that configuration is a fixed point of $\Phi$ and the trajectory stays there.
\end{proof}
Without growth every orbit of the deterministic rule ends on a fixed point, homogeneous or not, which is class-I or class-II behaviour in Wolfram's sense.

The next result concerns the interfaces between factions. Let $F(x)$ be the set of frontier cells and $\delta_i(x)$ the graph distance from $i$ to $F(x)$.

\begin{proposition}[Locality of interfaces]
\label{prop:local}
For either battle operator and any $g$, in one step of the local rule: (a) every cell with $\delta_i\ge3$ changes only through (G); (b) if the adjacency is symmetric, the states after the step of the cells with $\delta_i\le2$ do not depend on the troops of the cells with $\delta_i\ge3$: two configurations with the same owners whose troops agree on $\{\delta\le2\}$ have images that agree on $\{\delta\le2\}$ (given $\pi$ and, for dice, the random numbers).
\end{proposition}
\begin{proof}
Only frontier cells attack or are attacked, so owners change only on $F(x)$, and after (R) a cell can be a frontier cell only if $\delta_i\le1$. Senders in (F) are interior cells with a neighbour that is a frontier cell after (R), so they have $\delta_i\le2$; receivers have $\delta_i\le1$. Hence a cell with $\delta_i\ge3$ neither attacks, is attacked, sends nor receives, which is (a). For (b), the troops read in (D)--(F) when updating a cell with $\delta\le2$ are those of attackers and targets, which are frontier cells (a hostile neighbour of a frontier cell is itself a frontier cell when adjacency is symmetric), of senders, and of receivers, which are frontier cells after (R) and hence have $\delta\le1$ before it; all have $\delta\le2$. The remaining inputs are owners, which are the same in both configurations.
\end{proof}

Proposition~\ref{prop:local} says that the dynamics near an interface is locally determined, not that the interface stays in place. A one-dimensional CA requires two further conditions, neither of which follows from the rule. Suppose that over a time window the frontier cells of one interface stay within the strip of rows at distance at most $w$ from a fixed horizontal line on $\mathbb{T}_L$, and all other frontier cells stay at distance more than 10 from that strip. Then the rows within distance $w+2$ of the line evolve autonomously: each of their cells either changes only by growth (Proposition~\ref{prop:local}(a)) or lies within distance 2 of the interface, in which case its update depends on owners, which are constant outside the strip, and only on the troops of cells within distance 2 of the interface, which lie in these rows (Proposition~\ref{prop:local}(b) and Lemma~\ref{lem:radius}). On that window, and only on it, the front is a one-dimensional CA on $\mathbb{Z}_L$ whose cell is a column segment of height $2w+5$ and whose radius is at most 6 columns. Whether the conditions hold is an empirical question; E14 measures $w$.

\begin{hypothesis}[Complexity at intermediate attack margins]
\label{hyp:edge}
For $g\ge1$, the most complex behaviour (longest transients, intermediate damage spreading, most variable event entropy) occurs at intermediate attack margins, between fast consolidation at $\theta\approx0$ and saturation and freezing as $\theta\to\Nmax$. This is an ``edge of chaos'' hypothesis in the sense of \citet{langton1990}; it is tested in Section~\ref{sec:res:phase}.
\end{hypothesis}

\section{Experimental design}
\label{sec:design}

Defaults are $P=4$, $\Nmax=16$, $\rho=1.1718$, fortification on, synchronous updating, deterministic battles and quenched tie-breaking, with $o_i$ and $n_i$ initially independent and uniform. The experiments are listed in Table~\ref{tab:experiments}.

\begin{table}[t]
\centering\small
\caption{Experiments. $n_s$: runs per parameter combination; details in Section~\ref{sec:results}.}
\label{tab:experiments}
\begin{adjustbox}{max width=\textwidth}
\begin{tabular}{llllr}
\toprule
ID & Purpose & Board, rule & Factors & $n_s$ \\
\midrule
E1 & Battle kernel and its reduction & -- & $a,d\le30$ & exact \\
E2 & Class over $(g,\theta)$ & $\mathbb{T}_{32}$, det. & $g\le4$, $\theta\le12$, 2 tie-breaks & 16 \\
E3 & Finite-size scaling & $\mathbb{T}_L$, det. & 6 regimes, $L\le64$ & 12 \\
E4 & Damage spreading & $\mathbb{T}_{64}$ & 7 cases & 16 \\
E4b & Heterogeneity of damage & $\mathbb{T}_{64}$, det. & 3 regimes, 2 perturbations & 64 \\
E5 & Lattice vs.\ mean field & $\mathbb{T}_L$ & 4 variants, $L\le48$ & 16 \\
E5b & Random 4-regular graph & graph, $N\le2304$ & det./dice & 16 \\
E6 & Number of factions & $\mathbb{T}_{32}$ & $P\le8$, det./dice & 24 \\
E7 & 42-territory board & board, dice & $P\in\{3..6\}$, 2 schedules & 500 \\
E8 & Sensitivity analysis & $\mathbb{T}_{32}$, det. & 6 factors & 384 \\
E9 & Open boundaries & $32\times32$ & $P\in\{2,4\}$, det./dice & 24 \\
E10 & Moore neighbourhood & $\mathbb{T}_{32}$ & $g\le2$, $\theta\le12$ & 12 \\
E10b & Long E10 runs & $\mathbb{T}_{32}$, det. & $g=1$, $\theta\in\{8,12\}$, $T_{\max}=10^5$ & 12 \\
E11 & Stalemated fronts & $\mathbb{T}_L$, $P=2$ & $L\le48$, det./dice & 16 \\
E12 & Front Lyapunov speeds & $\mathbb{T}_L$, $P=2$, det. & $L\le96$ & 8--12 \\
E13 & Transition or crossover & $\mathbb{T}_L$, det. & $\theta\le8$, $L\le64$ & 48 \\
E13b & Extension of E13 & $\mathbb{T}_L$, det. & $\theta\in\{4,5,6\}$, $L\in\{96,128\}$ & 48 \\
E13c & Censored E13b runs continued & $\mathbb{T}_L$, det. & $\theta=4$, $T_{\max}=1.2\times10^5$ & 32 runs \\
E14 & Front geometry and bulk & $\mathbb{T}_{48}$, $P=2$, det. & 6 initial conditions & 12 \\
E15 & Fortification fraction & $\mathbb{T}_{32}$, det. & $f\in\{1/3,2/3\}$, $g\in\{1,2\}$, all $\theta$ & 16 \\
E16 & Dice over $(g,\theta)$ & $\mathbb{T}_{32}$, dice & $g\in\{1,2\}$, all $\theta$ & 16 \\
E17 & Continuous attack threshold & $\mathbb{T}_L$, det., random attacks at margin $\theta$ & $\theta_c\in[3.8,5.2]$, $L\le96$ & 24 \\
E18 & Interaction of oscillators & $\mathbb{T}_{32}$, det. & 104 long cycles of E2 & 40 tests \\
E19 & Larger and longer fronts & $\mathbb{T}_L$, $P=2$, det. & $L\in\{64,128\}$, $T\le5\times10^5$ & 8 \\
\bottomrule
\end{tabular}
\end{adjustbox}
\end{table}

\subsection{Measures}
\label{sec:measures}
For the deterministic rule the transient $\tau$ and period $p$ are exact (Proposition~\ref{prop:finite}). We store a 128-bit BLAKE2b digest of each configuration and stop at the first repeat; the probability of a digest collision over the whole study is below $10^{-25}$. Runs without a repeat within $T_{\max}$ are \emph{censored}. For cycles with $p>1$ we also record the period of each cell's own state sequence. The global period is always the least common multiple of these local periods; what distinguishes attractors is whether some cell carries the global period itself. If none does, the long period arises only from combining spatially fixed oscillators with shorter, different periods; if many cells do, the cycle is carried by structures that move across the lattice.

For damage spreading, the owner of one random cell is changed in a copy of the configuration and both copies are evolved, with shared random numbers for stochastic rules \cite{grassberger1995}. We record the fraction of cells with different owners, $H_o(t)$, or different state (owner or troops), $H_f(t)\ge H_o(t)$, and the damage radius. The two differ when a perturbation changes troop counts without changing the territorial outcome. Damage heals or stays confined in classes I and II, spreads linearly through the system in class III, and spreads irregularly in class IV \cite{wolfram1984,derrida1986}.

The compression ratio $C$ is the \texttt{zlib} size of the ownership space-time block for $t<512$, relative to a shuffle of the same symbols \cite{zenil2010}. For the event entropy, each cell at each step is assigned one of nine events (idle, grow, attack won or lost, defence held, conquered, fortification sent or received, held back); $\sigma_S$ is the standard deviation over $t<2000$ of the Shannon entropy of the event histogram, which, like Wuensche's input-entropy variance \cite{wuensche1999}, is large when quiescent and active periods alternate. We also record the activity, the fraction of cells that change owner in a step; the interface density $\rho_I$, the fraction of edges joining cells of different factions; the consolidation time $T_c$ (until a single owner remains); and the contagion ratio
$\kappa=\Pr[i\text{ conquered at }t{+}1\mid\exists j\in\mathcal{N}(i)\text{ conquered at }t]\,/\,\Pr[i\text{ conquered at }t{+}1]$,
computed over the first 400 steps, up to the last step with a conquest, with the neighbours of the board (lattice neighbours also for the mean-field control).

\subsection{Class criteria}
\label{sec:criteria}
Wolfram's classes describe the space-time behaviour of a rule from typical initial conditions, and they are asymptotic, whereas a simulation observes a finite transient, a finite-time response to perturbations and, if it ends, an attractor. We keep the three apart. Labelling runs by their attractor is an operational choice made for this study: a run labelled I or II has reached a homogeneous or a heterogeneous attractor, which says nothing about its transient and is not a claim about the Wolfram class of the rule as a whole. The transient and the response to perturbations are reported separately for that reason. Because I and II exhaust the possible attractors, the labels alone cannot reveal class-III or class-IV behaviour; that question is addressed by the structure of the attractors (whether a long cycle is carried by moving structures or by fixed local oscillators, Section~\ref{sec:res:artefact}), by the transients and damage, and by the fronts. Wolfram's classes were introduced for translation-invariant rules; with quenched priorities each run is a rule with fixed spatial disorder, and the labels describe the ensemble over $\pi$, not a single rule. The \emph{attractor class} is assigned only when the deterministic run reaches an attractor within $T_{\max}$: class I if it has a single owner (by Proposition~\ref{prop:fixed}(a), for $g\ge1$ this attractor is the saturated homogeneous fixed point), class II otherwise, whether a fixed point or a cycle of any period. Cycles with $p>L$ are reported separately as \emph{long cycles} because they distinguish the two tie-breaking rules (Section~\ref{sec:res:artefact}), but they are not a separate class. A censored run is \emph{unresolved} (U); it may have a long transient, a long period or neither. Unresolved runs are labelled by their finite-time sensitivity, U$_\mathrm{s}$ if the owner damage $H_o$, averaged over $t\in[351,400]$ after a perturbation at $t=0$, is at least 0.4 and U$_\mathrm{l}$ otherwise; the label describes a transient response and makes no claim about class III or IV. The transient length and $H_o$ are reported for all runs. A censored run is called a \emph{stalemate} if at least two factions survive, more than half of the rows or columns are owned by a single faction, and the mean activity over the last 500 steps exceeds $10^{-3}$. Each parameter combination takes the majority label of its runs. The thresholds were fixed after a pilot sweep; their effect is tested in Section~\ref{sec:res:phase}.

\subsection{Statistics and reproducibility}
Medians and means carry 95\% bootstrap intervals \cite{efron1993} and proportions Wilson intervals. Scaling exponents are least-squares fits on log--log data over the individual runs that reach an attractor, with bootstrap intervals; censored runs are reported separately. For comparison we also give the fit to the medians and a maximum-likelihood fit of all runs with a log-normal transient right-censored at $T_{\max}$ (a Tobit model). Groups are compared with Kruskal--Wallis, Mann--Whitney and Fisher tests, paired speeds with the Wilcoxon test, and trends with Spearman's $\rho_s$; Holm-adjusted $p$-values for all seventeen tests are in Appendix~\ref{app:tests}. Global sensitivity uses partial rank correlation coefficients on a Latin hypercube \cite{marino2008}.

On the 42-territory board (E7), territories are dealt at random and $P\in\{3,4,5,6\}$ factions start with the standard 35, 30, 25 or 20 units, one on each of their territories and the rest placed uniformly at random among them. We use $\Nmax=60$, $\theta=1$, dice and Risk-like reinforcement. A game ends with total conquest or after 600 rounds; a region is \emph{secured} when one faction holds it for three consecutive rounds.

The seed of every run is derived from the master seed 20261112 by \texttt{numpy}'s \texttt{SeedSequence}, keyed by a SHA-256 hash of the experiment name and parameters, so any single run can be regenerated. Unit tests check the dice probabilities, the board, Propositions~\ref{prop:balance}--\ref{prop:local}, Lemma~\ref{lem:radius} and the cycle detector. In addition, an independent implementation of the rule in Lean~4, with the battle operator as a parameter, machine-checks every formal result of Section~\ref{sec:model} for every board and state: Lemma~\ref{lem:radius} with its three witnesses and Propositions~\ref{prop:finite}--\ref{prop:local}, for deterministic battles and, where the statements include dice, for every realisation of the dice, including the dice operator of the code. Appendix~\ref{app:lean} lists the Lean declarations. Lean also checks the examples of Figure~\ref{fig:schematic} and Appendix~\ref{app:witness}, and the single-roll counts, $\rho$, the battle probabilities and the comparison of $\mathcal{R}_\rho$ with the exact kernel reported in Section~\ref{sec:res:dice}. The repository contains the per-run results of every experiment (JSON, 11\,MB), from which \texttt{analyse.py} regenerates all tables and figures without rerunning the simulations, and a machine-readable list of experiment configurations. The study runs in about nine hours on four cores with Python 3.13, NumPy 2.5, SciPy 1.18 and Matplotlib 3.11.

\section{Results}
\label{sec:results}

\subsection{Battle kernel (E1)}
\label{sec:res:dice}
\begin{figure}[t]
\centering
\includegraphics[width=0.48\textwidth]{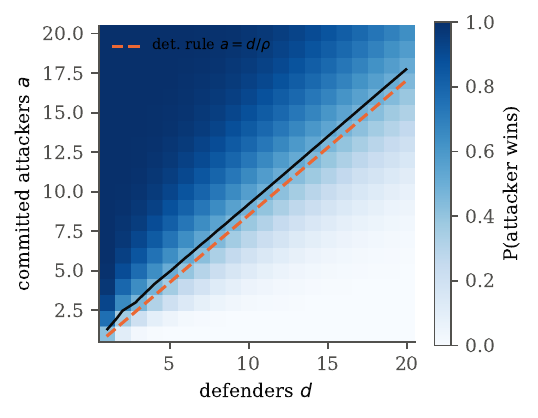}
\caption{Probability that $a$ committed units defeat $d$ defenders, with the 50\% contour (black) and the deterministic boundary $a=d/\rho$ (dashed).}
\label{fig:kernel}
\end{figure}

Enumeration reproduces Osborne's corrected single-roll probabilities \cite{osborne2003}: with three dice against two, the defender loses two units with probability $2890/7776$, the attacker two with $2275/7776$, and each one with $2611/7776$. The expected losses per roll, 0.9209 for the attacker and 1.0791 for the defender, give $\rho=1.1718$. The probability that $a$ committed units defeat $d$ defenders (Figure~\ref{fig:kernel}) is 0.417 for $(1,1)$, 0.506 for $(5,5)$ and 0.633 for $(20,20)$, so small equal battles favour the defender and large ones the attacker. $\mathcal{R}_\rho$ predicts the more probable outcome in 97.6\% of the grid $a,d\le30$, with a mean error of 0.81 units in the surviving attackers. All 22 errors lie next to the boundary $a=d/\rho$, at every battle size from $(1,1)$ to $(26,30)$; in each of them $\mathcal{R}_\rho$ awards the battle to the attacker although the exact probability of winning is 0.36--0.50. The board simulations use the exact kernel.

\subsection{Classification over the $(g,\theta)$ plane (E2)}
\label{sec:res:phase}
\begin{figure}[t]
\centering
\includegraphics[width=\textwidth]{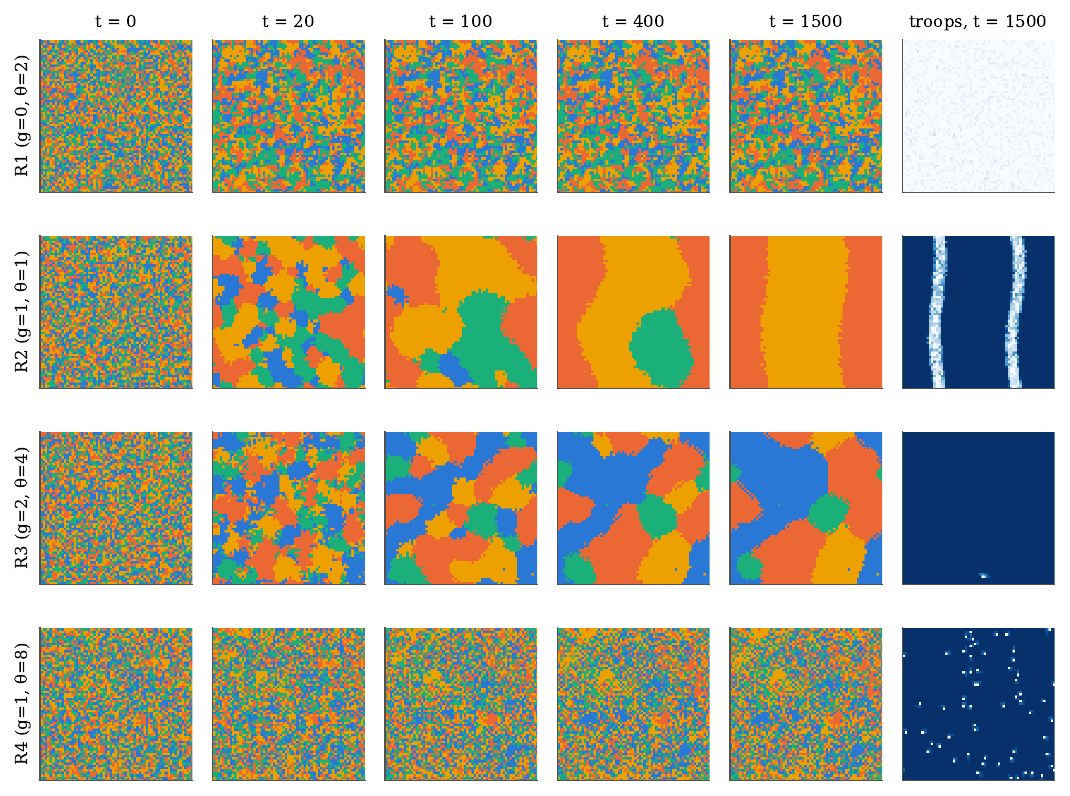}
\caption{Ownership at selected times (one colour per faction, $P=4$) and troops at $t=1500$ (last column; light $=$ few units) on $\mathbb{T}_{64}$. R1 ($g=0$) freezes at once. R2 ($g=1,\theta=1$) coarsens into compact domains with depleted borders; most R2 runs end with one owner, but this one ends in two straight stripes whose fronts are still active after $3\times10^4$ steps (a stalemate, Section~\ref{sec:res:robust}). R3 ($g=2,\theta=4$) coarsens into a partition that ends on a saturated fixed point or a cycle (31 and 28 of 60 runs in E3). R4 ($g=1,\theta=8$) keeps a fine mosaic with a few oscillating cells.}
\label{fig:snapshots}
\end{figure}

\begin{figure}[t]
\centering
\includegraphics[width=\textwidth]{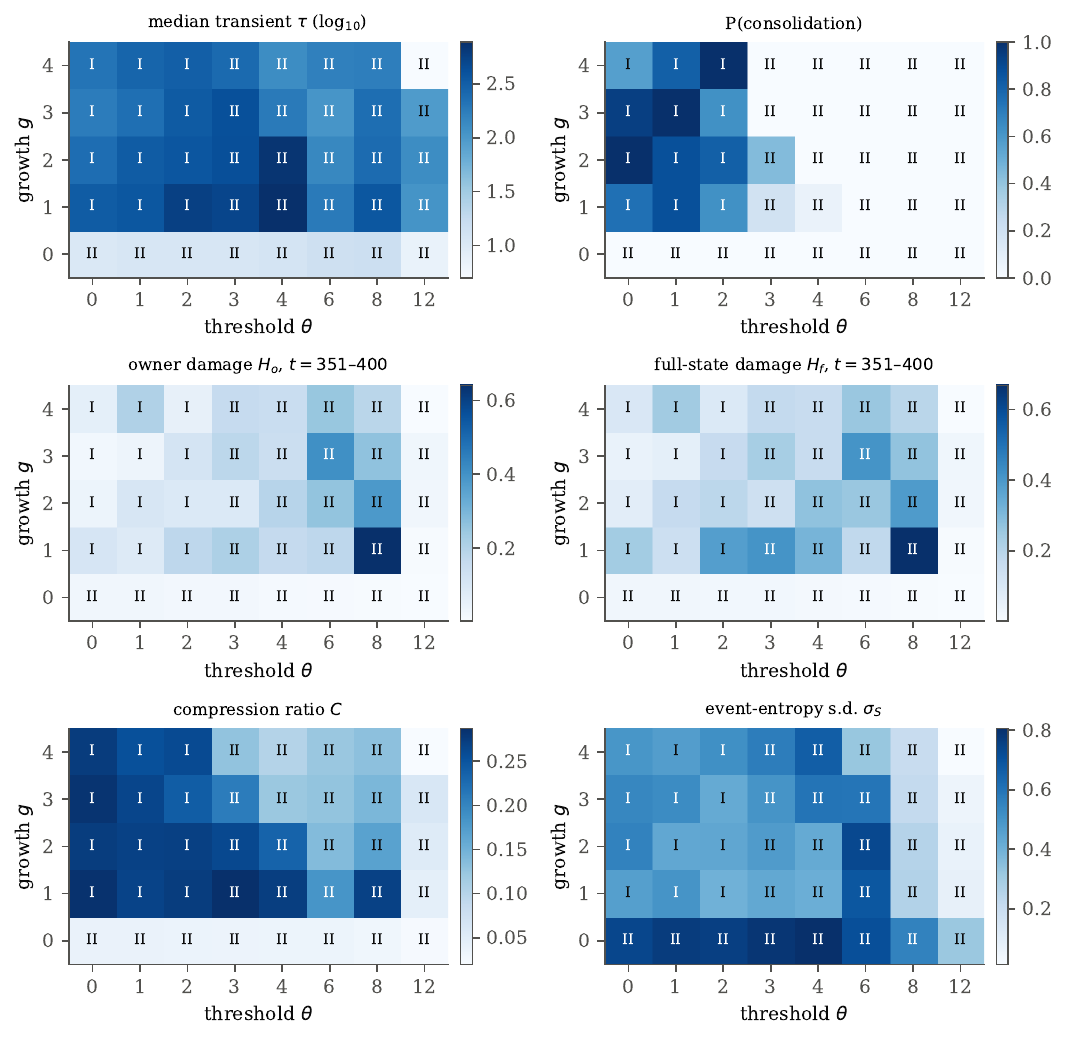}
\caption{Phase diagram on $\mathbb{T}_{32}$ (16 runs per cell). Colour: the quantity named above each panel (damage, $C$ and $\sigma_S$: mean over runs); labels: majority label (Section~\ref{sec:criteria}).}
\label{fig:phase}
\end{figure}

Figure~\ref{fig:snapshots} shows four typical regimes, Figure~\ref{fig:phase} the phase diagram from 640 runs ($L=32$, $T_{\max}=3\times10^4$), and Table~\ref{tab:perclass} the statistics of each label.

\begin{table}[t]
\centering\small
\caption{Medians by label (95\% bootstrap CI), E2.}
\label{tab:perclass}
\begin{tabular}{l r cccc}\toprule
class & $n$ & $\tilde\tau$ & $\tilde H_o$ & $\tilde C$ & $\tilde\sigma_S$\\ \midrule
I & 169 & 275 [254, 303] & 0.000 [0.000, 0.000] & 0.266 [0.263, 0.271] & 0.531 [0.510, 0.541]\\
II & 451 & 134 [123, 148] & 0.030 [0.020, 0.045] & 0.109 [0.101, 0.119] & 0.509 [0.449, 0.539]\\
U$_\mathrm{l}$ & 15 & $\ge$30000 (all censored) & 0.019 [0.000, 0.087] & 0.254 [0.240, 0.291] & 0.156 [0.096, 0.172]\\
U$_\mathrm{s}$ & 5 & $\ge$30000 (all censored) & 0.445 [0.416, 0.703] & 0.255 [0.239, 0.324] & 0.167 [0.102, 0.185]\\
\bottomrule\end{tabular}

\end{table}

For $g=0$ every run reaches a multi-owner fixed point within 6--14 steps (median), as Proposition~\ref{prop:g0} requires. These runs are more than a quarter of the class-II runs, which is why class II has shorter median transients than class I in Table~\ref{tab:perclass} (134 against 275 steps; 180 for class-II runs with $g\ge1$): reaching a single owner requires coarsening until one domain remains, whereas many multi-owner attractors are reached by freezing. For $g\ge1$ and $\theta\le2$, 56--100\% of runs end with one owner after median transients of 173--558 steps. Consolidation proceeds by domain-wall motion, and the walls appear in the troop field as depleted bands (Figure~\ref{fig:snapshots}, R2).

The probability of consolidation falls quickly with $\theta$: for $g=1$ it is 0.75, 0.88, 0.62, 0.19, 0.06 and 0 at $\theta=0,\dots,4$ and 6. Beyond this point runs freeze into saturated partitions; all 100 multi-owner fixed points with $g\ge1$ have $n_i\equiv\Nmax$, as in Proposition~\ref{prop:fixed}(b). The median transient peaks near the change (772 steps at $g=1,\theta=4$, against 330 at $\theta=0$ and 181 at $\theta=6$), consistent with Hypothesis~\ref{hyp:edge} for transients; the attractors on both sides are nevertheless fixed points. Section~\ref{sec:res:transition} analyses this change at larger sizes.

Of the 451 class-II runs, 111, almost all with $\theta\ge6$, end on long cycles, with periods from 48 to 27\,720 (Section~\ref{sec:res:artefact}). Only 20 runs (3.1\%) remain unresolved after $3\times10^4$ steps, 5 of them U$_\mathrm{s}$. Eighteen of these have $\theta\le2$, two or three factions and a low interface density, consistent with a few long domain walls; they are examined in E3 and E11. Owner damage is larger in class-II than in class-I runs (median 0.030 against 0; Mann--Whitney $p=4\times10^{-18}$), whereas compressibility shows no significant correlation with damage (Spearman $p=0.35$) and instead separates active coarsening ($C\approx0.26$) from frozen or oscillating mosaics ($C\approx0.03$--$0.13$). Neither is a reliable class indicator alone \cite{culik1988,zenil2010}. The unresolved runs have the compression ratio of coarsening runs ($C\approx0.25$) but the lowest entropy variability ($\sigma_S\approx0.16$), as expected for a few long domain walls that remain active while the rest of the lattice is quiescent.

The majority labels do not depend on the thresholds (Table~\ref{tab:sensitivity}). Moving the damage threshold between 0.2 and 0.6 relabels at most 4 of the 20 unresolved runs and no parameter cell. Shortening the window to 8000 or 3000 steps leaves 9 or 17 more runs unresolved but changes no majority label; only $T_{\max}=1000$ changes one. The long-cycle count depends on the cut-off (119, 111 and 81 runs for $p>L/2$, $L$ and $2L$), which is one reason not to treat long cycles as a class.

\subsection{Tie-breaking and long cycles}
\label{sec:res:artefact}
Repeating E2 with index-ordered tie-breaking, $\pi(i)=i$, barely changes the label counts (I/II/U$_\mathrm{s}$/U$_\mathrm{l}$: 169/451/5/15 quenched, 172/454/1/13 index), but the long cycles change character. With quenched priorities, the 111 long periods are highly composite (60, 84, 120, \dots, 27\,720) and only 6.3\% are multiples of $L$. Local periods were computed for the 104 cycles with $p\le8000$: in 94 of them no cell has the global period, the longest local period is in the median a third of the global one, the local periods take typically four distinct values, and a median of 5.6\% of cells oscillate. These attractors are combinations of spatially fixed local oscillators with different periods. E18 tests whether the oscillators interact. In each of the 104 cycles, the oscillating cells were grouped into clusters more than 6 cells apart (median 4 per cycle, range 1--8) and into contiguous oscillators (median 9, range 1--22). If the oscillators do not interact, setting each of them to an independent phase of the cycle gives another orbit, so that one step of $\Phi$ maps the shifted configuration to the shifted successor; this was checked exactly for 40 random choices of phases per cycle. The clusters passed every trial in all 104 cycles, and the contiguous oscillators in 95. In the other 9 cycles, all with $\theta\ge6$, 1--8 of the 40 trials failed, so some oscillators closer than the radius interact. In 81 of the 95 cycles no oscillator has the global period, which arises only from combining independent oscillators. Passing a sample of trials does not prove independence, but these long cycles are consistent with class II in the sense of separated periodic structures. With index order there are 59 long cycles, 58 of them with periods that are multiples of $L=32$; in 56 of the 57 analysed, some cells have the global period, and in the median case the only local periods are 1 and the global period. These are structures translating around the torus. They resemble the localised moving structures of class IV, but since they are absent with quenched priorities, they are produced by the asymmetry of the tie-breaking rule.

\subsection{Finite-size scaling (E3)}
\label{sec:res:scaling}

\begin{table}[t]
\centering\small
\caption{Finite-size scaling $\tau\propto L^z$ ($L=16$--64, 12 runs per $L$). Least-squares fit on the runs that reach an attractor (95\% bootstrap CI) and its $R^2$; fit to the medians (censored runs at $T_{\max}=3\times10^4$); maximum-likelihood fit with censoring and its residual standard deviation $\sigma$ (in $\ln\tau$). Stalemates: censored runs with two or more factions, mostly uniform rows or columns, and continuing combat. $^a$In R4 censoring is caused by long periods, not long transients, so the censored model does not apply, and the fit to the medians (2.61) reflects the growth of periods beyond $T_{\max}$, not of transients.}
\label{tab:scaling}
\begin{adjustbox}{max width=\textwidth}
\begin{tabular}{l c c c c c l c}\toprule
regime & $z$ (uncensored) & $R^2$ & $z$ (medians) & $z$ (censored ML) & $\sigma$ & censored runs per $L$ & stalemates\\ \midrule
R1 frozen (g=0,$\theta$=2) & 0.08 [0.05, 0.11] & 0.29 & 0.08 & 0.08 [0.05, 0.10] & 0.06 & [0, 0, 0, 0, 0] & 0\\
R2 coarsening (g=1,$\theta$=1) & 2.03 [1.81, 2.22] & 0.67 & 2.17 & 2.53 [1.93, 3.13] & 1.57 & [0, 1, 2, 1, 3] & 4\\
R3 partition (g=2,$\theta$=4) & 1.93 [1.70, 2.19] & 0.74 & 2.01 & 1.93 [1.70, 2.16] & 0.57 & [0, 0, 0, 0, 0] & 0\\
R4 high-threshold (g=1,$\theta$=8) & 0.40 [0.15, 0.64] & 0.19 & 2.61 & n/a$^a$ & -- & [0, 0, 0, 3, 7] & 0\\
R5 coarsening (g=2,$\theta$=0) & 2.13 [1.94, 2.33] & 0.79 & 2.09 & 1.94 [0.95, 2.93] & 1.84 & [2, 0, 2, 1, 2] & 4\\
R4-index (g=1,$\theta$=8) & 0.24 [-0.03, 0.53] & 0.05 & 0.13 & 0.24 [-0.06, 0.53] & 0.51 & [0, 0, 0, 0, 0] & 0\\
\bottomrule\end{tabular}

\end{adjustbox}
\end{table}

\begin{figure}[t]
\centering
\includegraphics[width=0.55\textwidth]{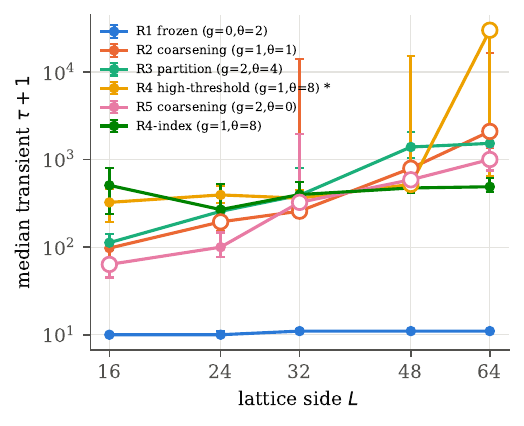}
\caption{Median transient against $L$ (95\% bootstrap CI), with censored runs counted at $T_{\max}$. Open circles: sizes with at least one censored run. $^*$For R4 the censoring at $L\ge48$ is caused by periods longer than $T_{\max}$, so the rise of this curve reflects period growth, not longer transients.}
\label{fig:scaling}
\end{figure}

The frozen regime R1 has an almost size-independent transient (median 9--10 steps, $z=0.08$ [0.05, 0.11]; Table~\ref{tab:scaling}, Figure~\ref{fig:scaling}). For the coarsening regimes R2, R5 and R3 we report two strata separately. Among runs that reach an attractor, the transient grows as $L^{2.03}$, $L^{2.13}$ and $L^{1.93}$ (fits to the medians: 2.17, 2.09, 2.01), consistent with the diffusive $\tau\propto L^2$ of curvature-driven domain coarsening \cite{bray1994}. The censored runs form a second stratum: of the 14 in R2 and R5, 8 are stalemated fronts (Section~\ref{sec:res:robust}), which need not end at all, and their share grows with $L$ (Table~\ref{tab:scaling}). A single log-normal model with censoring, which mixes the two strata, gives 2.53 [1.93, 3.13], 1.94 [0.95, 2.93] and 1.93 [1.70, 2.16]; we report it for comparison only, since a competing-risks description (consolidation versus stalemate) would be more appropriate and needs more runs than we have. The run-to-run scatter is large ($R^2=0.67$--0.79 in $\ln\tau$), and five sizes with 12 runs each provide a weak consistency check with $L^2$, not an estimate of an asymptotic exponent. The per-run data are in the repository. In R4 the period rather than the transient grows with size, from at most 720 at $L=16$ to 21\,840 at $L=32$, with ten of 24 runs censored for $L\ge48$. This is consistent with the oscillator decomposition, since the least common multiple of $k$ small integers grows roughly exponentially with $k$. With index order, R4 transients do not grow ($z=0.24$) and all 15 long periods are multiples of the respective $L$ ($208=13\cdot16$, $864=27\cdot32$, $3392=53\cdot64$), confirming the tie-breaking artefact at all sizes.

\subsection{Dominance--fragmentation crossover (E13--E13c)}
\label{sec:res:transition}
\begin{figure}[t]
\centering
\includegraphics[width=\textwidth]{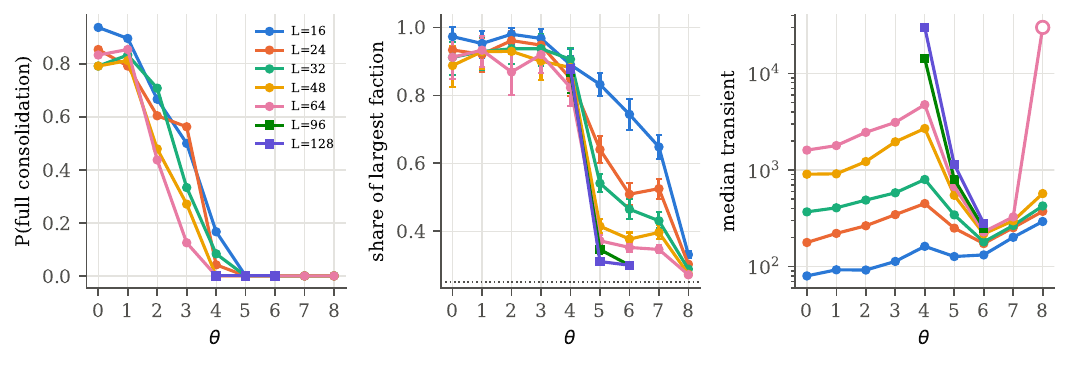}
\caption{E13--E13c (censored runs at $\theta=4$, $L\ge96$ continued to $1.2\times10^5$ steps): probability of full consolidation, share of the largest faction (dotted: $1/P$; bars: 95\% bootstrap CI) and median transient (censored runs at $T_{\max}$) against $\theta$. Open symbols: points where at least half of the runs are censored; the jump at $\theta=8$, $L=64$ is caused by periods longer than $T_{\max}$, not by long transients, as in R4 of E3.}
\label{fig:transition}
\end{figure}

To test whether the change found in E2 is a phase transition, the deterministic rule ($P=4$, $g=1$) was run for $\theta=0,\dots,8$ at $L=16$ to 64, and for $\theta\in\{4,5,6\}$ also at $L=96$ and 128 (E13b), with 48 runs per point (Table~\ref{tab:transition}, Figure~\ref{fig:transition}).

\begin{table}[t]
\centering\small
\caption{E13 and E13b: share of the largest faction and, in brackets, probability of full consolidation (48 runs per entry; --: not run; at $\theta=4$, $L\ge96$, runs censored at $3\times10^4$ are replaced by their continuation to $1.2\times10^5$ steps, E13c). $\theta_{50}$ and $w$: midpoint and scale of a logistic fit $1/(1+e^{(\theta-\theta_{50})/w})$ of the consolidation probability, whose 10--90\% width is $4.4\,w$ (95\% CI).}
\label{tab:transition}
\begin{adjustbox}{max width=\textwidth}
\begin{tabular}{r ccccccccc c c}\toprule
 & \multicolumn{9}{c}{share of largest faction (P(full consolidation))} & & \\
$L$ & $\theta=0$ & $\theta=1$ & $\theta=2$ & $\theta=3$ & $\theta=4$ & $\theta=5$ & $\theta=6$ & $\theta=7$ & $\theta=8$ & $\theta_{50}$ & scale $w$\\ \midrule
16 & 0.97 (0.94) & 0.95 (0.90) & 0.98 (0.67) & 0.97 (0.50) & 0.89 (0.17) & 0.83 (0.00) & 0.74 (0.00) & 0.65 (0.00) & 0.33 (0.00) & 2.66 [2.38, 2.91] & 0.77 [0.65, 0.98]\\
24 & 0.93 (0.85) & 0.92 (0.79) & 0.96 (0.60) & 0.95 (0.56) & 0.85 (0.04) & 0.64 (0.00) & 0.51 (0.00) & 0.53 (0.00) & 0.30 (0.00) & 2.32 [2.07, 2.59] & 0.89 [0.75, 1.06]\\
32 & 0.91 (0.79) & 0.93 (0.83) & 0.94 (0.71) & 0.94 (0.33) & 0.91 (0.08) & 0.54 (0.00) & 0.46 (0.00) & 0.43 (0.00) & 0.29 (0.00) & 2.21 [1.94, 2.44] & 0.88 [0.71, 1.06]\\
48 & 0.89 (0.79) & 0.93 (0.81) & 0.93 (0.48) & 0.90 (0.27) & 0.88 (0.00) & 0.41 (0.00) & 0.38 (0.00) & 0.40 (0.00) & 0.27 (0.00) & 1.81 [1.52, 2.05] & 0.81 [0.66, 0.98]\\
64 & 0.91 (0.83) & 0.93 (0.85) & 0.87 (0.44) & 0.92 (0.12) & 0.82 (0.00) & 0.37 (0.00) & 0.35 (0.00) & 0.35 (0.00) & 0.27 (0.00) & 1.73 [1.46, 1.98] & 0.66 [0.52, 0.81]\\
96 & -- & -- & -- & -- & 0.87 (0.00) & 0.35 (0.00) & 0.30 (0.00) & -- & -- & -- & --\\
128 & -- & -- & -- & -- & 0.88 (0.00) & 0.31 (0.00) & 0.30 (0.00) & -- & -- & -- & --\\
\bottomrule\end{tabular}

\end{adjustbox}
\end{table}

The probability of full consolidation does not sharpen with size. Its logistic midpoint drifts from $\theta_{50}=2.66$ [2.38, 2.91] at $L=16$ to $1.73$ [1.46, 1.98] at $L=64$, while the scale of the fit stays nearly constant ($w\propto L^{-0.11}$), and at $L\ge96$ no run consolidates for $\theta\ge4$. The reason is the saturation of minority domains: of the 117 runs at $L=64$ with $\theta=2$--4 that do not consolidate, 91 end on multi-owner fixed points, in which the largest faction holds a median 99.9\% of the cells; all such fixed points observed in E2 are saturated, and a saturated cell cannot be attacked when $\theta\ge1$, since condition~\eqref{eq:attack} would require $n_i\ge\Nmax+\theta$. Larger lattices contain more minority domains, so the chance that at least one survives grows with $L$ even when one faction holds most of the lattice.

The share of the largest faction behaves differently. At $\theta=5$ it falls from 0.83 [0.80, 0.86] at $L=16$ to 0.37 [0.36, 0.39] at $L=64$ and then levels off (0.35 and 0.31 at $L=96$ and 128, against a floor of $1/P=0.25$). At $\theta=4$, 26 of the 48 runs at $L=128$ (and 6 at $L=96$) were still coarsening at $T_{\max}=3\times10^4$, which biases the share downwards (0.77 at $L=128$). Continuing these runs to $1.2\times10^5$ steps (E13c), 23 of the 26 at $L=128$ reach a multi-owner attractor and none consolidates. Over the 45 runs at $L=128$ that reached an attractor the share is 0.90 [0.86, 0.93]; counting the 3 unresolved runs at their state after $1.2\times10^5$ steps (shares 0.46--0.76), it is 0.88 [0.84, 0.91]. At $L=96$ the corresponding values are 0.92 [0.88, 0.95] over 43 resolved runs and 0.87 [0.81, 0.91] including the 5 unresolved ones (shares 0.36--0.49, consistent with stalemated fronts). Over $L=16$--128 the share at $\theta=4$ thus stays at 0.82--0.91 while that at $\theta=5$ falls towards $1/P$: the step in the largest faction's share lies between $\theta=4$ and 5 at every size and becomes steeper with $L$. The drop in the share between $\theta=4$ and 5 is 0.06, 0.21, 0.36, 0.47, 0.45, 0.52 and 0.57 for $L=16$, 24, 32, 48, 64, 96 and 128; it grows quickly up to $L=48$ and slowly afterwards, and with an integer $\theta$ the width of the step cannot be resolved below one unit. The median transient (censored runs at $T_{\max}$) grows as $L^{2.1}$--$L^{2.5}$ for $\theta\le4$ but only as $L^{1.2}$ at $\theta=5$.

\begin{figure}[t]
\centering
\includegraphics[width=\textwidth]{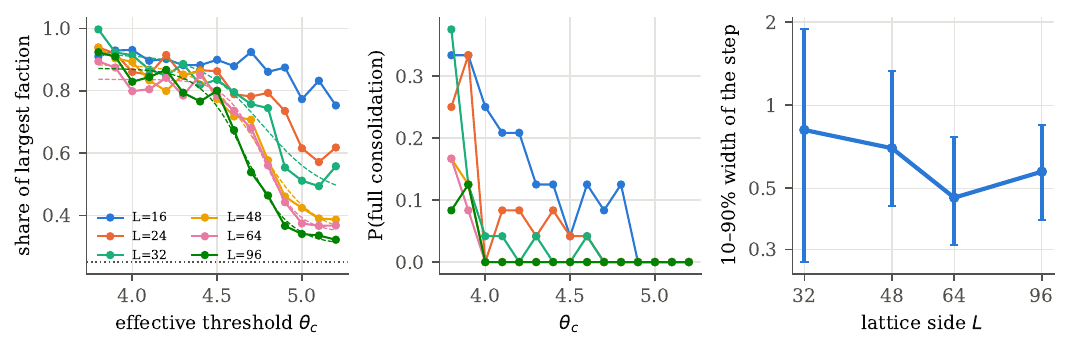}
\caption{E17: share of the largest faction (left; dashed: logistic fits; dotted: $1/P$) and probability of full consolidation (centre) against the continuous threshold $\theta_c$, and 10--90\% width of the step against $L$ (right, 95\% CI).}
\label{fig:contthreshold}
\end{figure}

E17 resolves the step with a continuous threshold. An attack whose margin is exactly $\theta$ is made only with probability $1-f$, so that the effective threshold $\theta_c=\theta+f$ interpolates between integers; the rule is then stochastic, and at integer $\theta_c$ it is the deterministic rule of E13, whose shares it reproduces (at $L=32$, 0.91 and 0.51 at $\theta_c=4$ and 5, against 0.91 and 0.54). Scanning $\theta_c$ from 3.8 to 5.2 in steps of 0.1, with 24 runs per point and $L=16$--96 (Table~\ref{tab:contthreshold}, Figure~\ref{fig:contthreshold}), places the drop of the share at $\theta^*=4.67$--4.77 for $L=32$--96, with a 10--90\% width of 0.46--0.81, below one unit of $\theta$ (for $L\ge64$ also at the upper end of the 95\% CI); for $L=16$ and 24 the drop extends beyond $\theta_c=5.2$, as in E13. Runs censored at $3\times10^4$ steps (at most 33 of 360 per size) enter with their final share. Between $L=32$ and 96 the midpoint moves by at most 0.1, but the width does not decrease significantly (width $\propto L^{-0.39}$, 95\% CI of the exponent $[-1.25, 0.59]$).

\begin{table}[t]
\centering\small
\caption{E17: continuous threshold $\theta_c$. Logistic fit of the share of the largest faction against $\theta_c\in[3.8,5.2]$ (24 runs per point, 95\% bootstrap CI); midpoint $\theta^*$, 10--90\% width of the drop and plateaus. For $L<32$ the drop lies beyond the scanned range.}
\label{tab:contthreshold}
\begin{tabular}{r c c c c r}\toprule
$L$ & $\theta^*$ (95\% CI) & 10--90\% width (95\% CI) & upper plateau & lower plateau & censored\\ \midrule
32 & 4.77 [4.72, 4.97] & 0.81 [0.27, 1.88] & 0.92 & 0.46 & 3/360\\
48 & 4.76 [4.71, 4.84] & 0.70 [0.43, 1.33] & 0.88 & 0.34 & 21/360\\
64 & 4.77 [4.73, 4.80] & 0.46 [0.31, 0.76] & 0.84 & 0.34 & 24/360\\
96 & 4.67 [4.63, 4.70] & 0.57 [0.38, 0.85] & 0.87 & 0.30 & 33/360\\
\bottomrule\end{tabular}

\end{table}

In short, the probability of full consolidation is a size-dependent crossover whose midpoint drifts with $L$, governed by the statistics of saturated minority domains, whereas the largest faction's share drops in a step that sharpens with $L$ up to $L=48$ without moving, and that, resolved with a continuous threshold, is 0.5--0.8 wide (10--90\%) at $\theta^*\approx4.7$ for $L=32$--96. The widths up to $L=96$ do not show the narrowing expected at a phase transition, nor do they exclude it, and we make no claim of a phase transition.

\subsection{Damage spreading (E4)}
\label{sec:res:damage}

\begin{figure}[t]
\centering
\includegraphics[width=\textwidth]{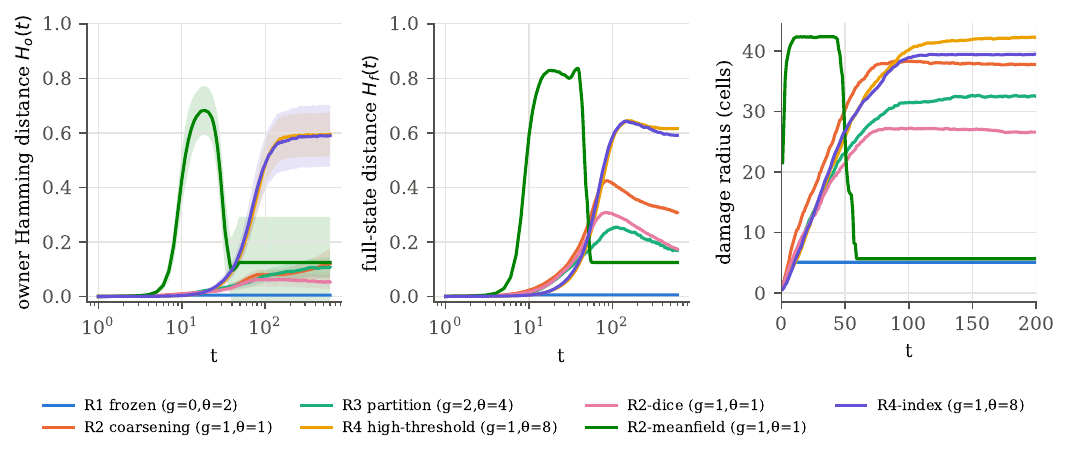}
\caption{Damage after changing one owner at $t=0$ on $\mathbb{T}_{64}$ (16 pairs; 95\% CI on the left). Left: fraction of cells with a different owner, $H_o$. Centre: fraction with a different owner or troop count, $H_f$; the gap to $H_o$ is damage to force levels that has not changed territory. Right: damage radius (not physical for the mean-field control).}
\label{fig:damage}
\end{figure}

\begin{table}[t]
\centering\small
\caption{Damage spreading (mean, 95\% CI) and initial speed $v$ of the damage front, the slope of the mean damage radius over $t\in[5,60]$.}
\label{tab:damage}
\begin{adjustbox}{max width=\textwidth}
\begin{tabular}{l ccc c}\toprule
case & $H_o(100)$ & $H_o(600)$ & $H_f(600)$ & front speed $v$ (cells/step)\\ \midrule
R1 frozen (g=0,$\theta$=2) & 0.005 [0.004, 0.007] & 0.005 [0.004, 0.007] & 0.006 [0.004, 0.008] & 0.003\\
R2 coarsening (g=1,$\theta$=1) & 0.080 [0.061, 0.096] & 0.121 [0.070, 0.174] & 0.308 [0.230, 0.381] & 0.532\\
R3 partition (g=2,$\theta$=4) & 0.074 [0.051, 0.096] & 0.107 [0.071, 0.143] & 0.169 [0.113, 0.224] & 0.417\\
R4 high-threshold (g=1,$\theta$=8) & 0.482 [0.397, 0.543] & 0.594 [0.507, 0.648] & 0.616 [0.525, 0.671] & 0.524\\
R2-dice (g=1,$\theta$=1) & 0.062 [0.037, 0.086] & 0.054 [0.028, 0.080] & 0.174 [0.102, 0.248] & 0.374\\
R2-meanfield (g=1,$\theta$=1) & 0.125 [0.000, 0.312] & 0.125 [0.000, 0.312] & 0.125 [0.000, 0.312] & --\\
R4-index (g=1,$\theta$=8) & 0.486 [0.392, 0.569] & 0.590 [0.468, 0.678] & 0.591 [0.470, 0.681] & 0.535\\
\bottomrule\end{tabular}

\end{adjustbox}
\end{table}

In R1 damage stays within about 5 cells ($H_o<0.007$). With growth it spreads at a roughly constant 0.37--0.54 cells per step, at most 11\% of the bound of Lemma~\ref{lem:radius} (5 cells per step for the deterministic cases, all with $\theta\ge1$, and 6 with dice; Figure~\ref{fig:damage}, Table~\ref{tab:damage}). The owner and full-state distances differ: in the coarsening regimes the perturbation changes the owner of only 5--12\% of cells by $t=600$ but the troops of 17--31\%, so most of the damage alters force levels without changing the territorial outcome. In R4 owner damage reaches $H_o=0.59$ [0.51, 0.65], although the attractor is of class II, i.e.\ sensitivity to initial conditions is confined to the transient. In the mean-field control the outcome is binary: 14 of 16 pairs end with the same winner, the other two with different winners.

The means hide a broad distribution (E4b, Table~\ref{tab:damagehet}). In R2 and R3, 19--39\% of perturbations heal completely, and in R2 the 90th percentile of $H_o$ is 3--4 times the median. In R4 an owner change almost always spreads ($H_o\approx0.65$), whereas a one-unit change heals in 61\% of pairs and otherwise spreads to the same level, so the response is bimodal. The final damage depends little on the troops or the number of hostile neighbours of the perturbed cell ($|\rho_s|\le0.34$).

\subsection{Spatial structure versus mean field and random graphs (E5, E5b)}
\label{sec:res:spatial}

\begin{figure}[t]
\centering
\includegraphics[width=\textwidth]{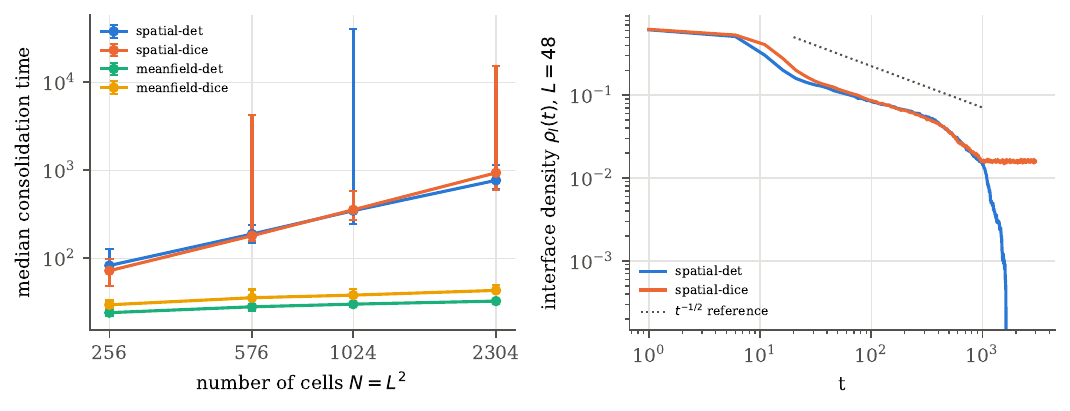}
\caption{Left: median time to consolidation against $N$ ($g=1,\theta=1$, 16 runs per point). Right: interface density on $\mathbb{T}_{48}$ with a $t^{-1/2}$ reference.}
\label{fig:spatial}
\end{figure}

\begin{table}[t]
\centering\small
\caption{Lattice, mean-field control and random 4-regular graph (E5b): scaling of consolidation time with $N$ (runs censored at $T_{\max}=4\times10^4$), interface decay exponent $\alpha$ (range over fitting windows), and contagion ratio $\kappa$ ($N=1024$; neighbours are lattice neighbours for the mean-field control and graph neighbours for the random graph).}
\label{tab:spatial}
\begin{adjustbox}{max width=\textwidth}
\begin{tabular}{l c c c c c}\toprule
variant & $z_N$ ($T_{c}\propto N^{z_N}$) & 95\% CI & censored & coarsening $\alpha$ (range over windows) & contagion $\kappa$ ($N=1024$)\\ \midrule
spatial-det & 1.09 & [0.95, 1.22] & 10/64 & 0.67 (0.39--0.74) & 5.22 [4.86, 5.59]\\
spatial-dice & 1.03 & [0.57, 1.44] & 8/64 & 0.67 (0.46--0.72) & 4.92 [4.66, 5.17]\\
meanfield-det & 0.14 & [0.10, 0.19] & 0/64 & -- & 1.11 [1.10, 1.13]\\
meanfield-dice & 0.17 & [0.12, 0.22] & 0/64 & -- & 1.15 [1.12, 1.17]\\
random 4-regular, det & 0.18 & [0.10, 0.26] & 11/64 & -- & 1.50 [1.48, 1.53]\\
random 4-regular, dice & 0.17 & [0.10, 0.23] & 5/64 & -- & 1.50 [1.47, 1.53]\\
\bottomrule\end{tabular}

\end{adjustbox}
\end{table}

In the mean-field control consolidation takes 24--43 steps and barely depends on size ($z_N\approx0.15$). On the lattice it scales as $N^{1.03}$ (dice) and $N^{1.09}$ (deterministic), i.e.\ as $L^2$, in agreement with E3 (Table~\ref{tab:spatial}, Figure~\ref{fig:spatial}). At $N=2304$ the lattice conflict lasts more than twenty times longer (934 against 43 steps with dice; at $L=32$, Mann--Whitney $p=8\times10^{-7}$). Conquests cluster in space: a cell next to a fresh conquest is $\kappa=4.9$--5.2 times more likely than average to be conquered next, against $\kappa\approx1.1$ in the control. The interface density decays approximately as a power law before the last domains vanish. Its exponent depends on the fitting window (0.39--0.74) and is not determined, but the decay is faster than the logarithmic decay of the voter model \cite{dornic2001}, consistent with curvature-driven coarsening \cite{bray1994}.

The mean-field control changes two things at once: it removes geometry and it replaces local by global mixing. E5b separates them by running the unchanged local rule on quenched random 4-regular graphs, which have the degree of the lattice and only local interactions but no low-dimensional geometry. Consolidation there is as fast as in the mean-field control (median 25--43 steps for $N=256$--2304, $z_N=0.18$ [0.10, 0.26] with deterministic battles and 0.17 with dice), and conquests are only weakly clustered ($\kappa=1.5$). The slow, spatially correlated conflict on the lattice is a consequence of its geometry rather than of local interaction as such. On the random graph, 16 of 128 runs end as frozen two-faction partitions.

\subsection{Number of factions (E6)}
\label{sec:res:players}
\begin{figure}[t]
\centering
\includegraphics[width=0.5\textwidth]{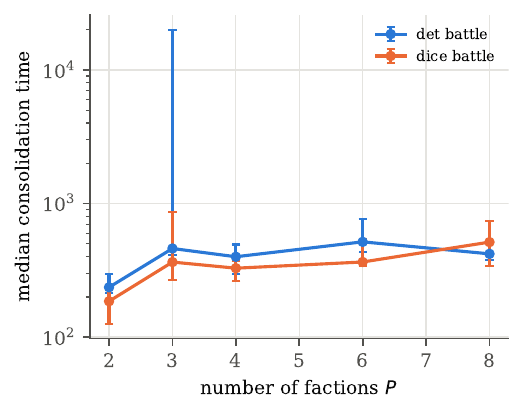}
\caption{Median time to consolidation against the number of factions ($\mathbb{T}_{32}$, $g=1$, $\theta=1$, 24 runs per point; censored runs counted at $T_{\max}=2\times10^4$; 95\% bootstrap CI).}
\label{fig:players}
\end{figure}

With $g=1$ and $\theta=1$, the median consolidation time with dice (censored runs counted at $T_{\max}=2\times10^4$) tends to increase with $P$, from 186 steps at $P=2$ to 514 at $P=8$ (Kruskal--Wallis $p=0.002$, Spearman Holm $p=0.007$; Figure~\ref{fig:players}). With deterministic battles the medians differ across $P$ (Kruskal--Wallis Holm $p=0.047$) but show no monotone trend (Spearman Holm $p=0.29$), because 7 and 11 of 24 runs at $P=2$ and 3 do not consolidate within $T_{\max}$. In an independent replicate of this setting with $P=2$ (E9), all 7 unconsolidated runs are stalemates (Section~\ref{sec:res:robust}).

\subsection{The 42-territory board (E7)}
\label{sec:res:board}
E7 runs the rule on the irregular board of the game, with dice and Risk-like reinforcement, as a qualitative check rather than a model of play: all factions follow one fixed heuristic, and the reinforcement makes the board model non-local.

\begin{table}[t]
\centering\small
\caption{Game length in rounds, 500 games per configuration.}
\label{tab:board}
\begin{tabular}{lrl r l}\toprule
schedule / $P$ / topology & finished & mean length [95\% CI] & median & IQR\\ \midrule
turn / 3 / spatial & 500/500 & 17.1 [16.7, 17.6] & 16 & 14--20\\
turn / 4 / spatial & 500/500 & 23.0 [22.4, 23.7] & 22 & 18--26\\
turn / 5 / spatial & 500/500 & 29.6 [28.7, 30.5] & 28 & 22--34\\
turn / 6 / spatial & 500/500 & 35.9 [34.8, 37.2] & 33 & 27--42\\
sync / 3 / spatial & 500/500 & 17.9 [17.5, 18.4] & 17 & 14--20\\
sync / 4 / spatial & 500/500 & 24.7 [24.1, 25.4] & 23 & 19--29\\
sync / 5 / spatial & 500/500 & 33.7 [32.6, 34.8] & 30 & 25--39\\
sync / 6 / spatial & 500/500 & 46.4 [44.8, 48.2] & 42 & 32--56\\
sync / 4 / meanfield & 500/500 & 85.4 [82.8, 88.1] & 82 & 65--99\\
turn / 4 / meanfield & 500/500 & 53.2 [51.2, 55.3] & 49 & 37--64\\
\bottomrule\end{tabular}

\end{table}

\begin{figure}[t]
\centering
\includegraphics[width=\textwidth]{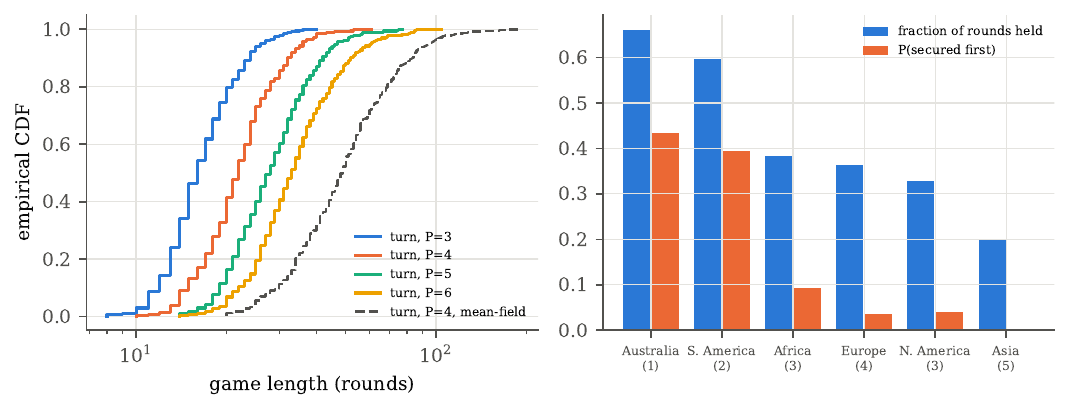}
\caption{Left: distribution of game length (turn schedule; dashed: mean-field control at $P=4$). Right: fraction of rounds in which each region is held by one faction, and probability that it is the first region secured ($P=4$, turn schedule); the number of entry cells is given in brackets.}
\label{fig:board}
\end{figure}

All 5000 games ended in total conquest, with median lengths of 16, 22, 28 and 33 rounds for $P=3$ to 6 (Table~\ref{tab:board}). Synchronous updating lengthens games, slightly at $P=4$ (mean 24.7 against 23.0 rounds, $p=6\times10^{-4}$) and more for six factions (46.4 against 35.9), without qualitative change \cite{schonfisch1999}. Small regions are held more often: Australia (one entry cell) is held in 66\% of rounds and is the first region secured in 43\% of games, Asia (five entry cells) in 20\% of rounds and almost never first (Figure~\ref{fig:board}). Across the six regions, the fraction of rounds held decreases with region size (Spearman $\rho_s=-0.99$, Holm $p=0.004$); its correlation with the number of entry cells is not significant after correction (Holm $p=0.10$). With six regions these are descriptive statements. On the board, unlike the torus, the mean-field control lengthens games (median 49 against 22 rounds, $p<10^{-130}$). A likely reason, which we have not tested, is a mismatch specific to the board model: reinforcement is placed on cells that are frontier cells on the board, while in the control attacks are made against random partners, so the reinforced cells are not the ones that fight. On the torus growth is uniform and there is no such mismatch.

\subsection{Sensitivity analysis (E8)}
\label{sec:res:sens}

\begin{figure}[t]
\centering
\includegraphics[width=\textwidth]{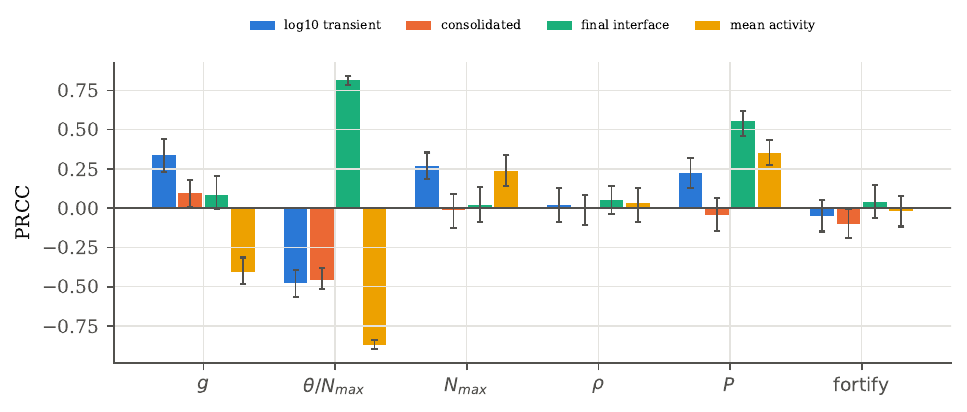}
\caption{Partial rank correlation coefficients (95\% CI) from 384 Latin-hypercube runs on $\mathbb{T}_{32}$ with the deterministic rule. Factor ranges: $g\in\{0,\dots,4\}$, $\Nmax\in\{8,\dots,32\}$, $\theta\in\{0,\dots,\Nmax-2\}$, $\rho\in[0.8,1.6]$, $P\in\{2,\dots,8\}$, fortification on or off; transients are censored at $6000$ steps.}
\label{fig:sens}
\end{figure}

The relative attack margin $\theta/\Nmax$ dominates every output (Figure~\ref{fig:sens}, Table~\ref{tab:prcc}): activity ($-0.87$), final interface density ($+0.81$), transient length ($-0.48$) and consolidation ($-0.46$). Growth comes second, lengthening transients ($+0.34$) and reducing activity ($-0.41$), and more factions increase fragmentation ($+0.55$). Fortification has a marginal effect on consolidation only, and the exchange ratio $\rho\in[0.8,1.6]$ has none. These conclusions hold within the sampled parameter ranges.

\subsection{Robustness and stalemated fronts (E9--E12, E14)}
\label{sec:res:robust}

\begin{table}[t]
\centering\small
\caption{E9: periodic and open $32\times32$ lattices ($g=1,\theta=1$, 24 runs per row). Active fronts: runs not consolidated after $2\times10^4$ steps in which more than $10^{-3}$ of the cells change owner per step over the last 500 steps; frozen: runs not consolidated in which fewer do.}
\label{tab:boundaries}
\begin{adjustbox}{max width=\textwidth}
\begin{tabular}{llr rrr c r}\toprule
boundary & battle & $P$ & not consolidated & active fronts & frozen & P(active front) [95\% CI] & median $T_{c}$\\ \midrule
periodic & det & 2 & 7/24 & 7 & 0 & 0.29 [0.15, 0.49] & 185\\
periodic & det & 4 & 5/24 & 3 & 2 & 0.12 [0.04, 0.31] & 399\\
periodic & dice & 2 & 7/24 & 7 & 0 & 0.29 [0.15, 0.49] & 192\\
periodic & dice & 4 & 4/24 & 3 & 1 & 0.12 [0.04, 0.31] & 316\\
open & det & 2 & 10/24 & 6 & 4 & 0.25 [0.12, 0.45] & 451\\
open & det & 4 & 9/24 & 3 & 6 & 0.12 [0.04, 0.31] & 695\\
open & dice & 2 & 17/24 & 13 & 4 & 0.54 [0.35, 0.72] & 411\\
open & dice & 4 & 5/24 & 4 & 1 & 0.17 [0.07, 0.36] & 951\\
\bottomrule\end{tabular}

\end{adjustbox}
\end{table}

\paragraph{Boundaries (E9).} With open boundaries, 26 of 96 runs remain unconsolidated with continuing combat, against 20 of 96 on the torus (Fisher $p=0.40$; Table~\ref{tab:boundaries}); in 44 of these 46 runs more than half of the rows or columns belong to a single faction, i.e.\ they are stalemated fronts. Stalemates do not require periodic boundaries. Open boundaries also leave more unconsolidated runs without owner changes (15 against 3) and slow consolidation (median 695 against 399 steps for $P=4$ with deterministic battles).

\paragraph{Neighbourhood (E10).} With the Moore neighbourhood the majority label agrees with the von Neumann result in 22 of 24 cells (Table~\ref{tab:moore}); the two exceptions, at $\theta=3$, consolidate more often with eight neighbours (75\% and 92\% of runs) and so change from II to I; as on the von Neumann lattice, multi-owner fixed points are saturated. Continuing the 24 runs at $g=1,\theta\in\{8,12\}$ to $10^5$ steps (E10b), 23 reach cycles after transients of up to 22\,985 steps. Nine of these transients exceed 8000 steps, so shorter observation windows would have left these runs unresolved.

\paragraph{Fortification fraction (E15).} Fortification sends half of the surplus, $\lfloor(n_i-n_{k^*})/2\rfloor$. Repeating the E2 runs for $g\in\{1,2\}$ with the same seeds but sending $\lfloor f(n_i-n_{k^*})\rfloor$ with $f=1/3$ or $2/3$ (Table~\ref{tab:fort}) changes the majority label of 2 ($f=1/3$) and 0 ($f=2/3$) of the 16 parameter cells; the two changes, at $(g,\theta)=(1,2)$ and $(2,3)$, lie where the consolidation probability is close to one half. The overall consolidation probability stays at 0.34--0.35, long cycles remain (52--58 runs), and every multi-owner fixed point is saturated for all three fractions. The freezing mechanism does not depend on the fraction.

\begin{figure}[t]
\centering
\includegraphics[width=0.55\textwidth]{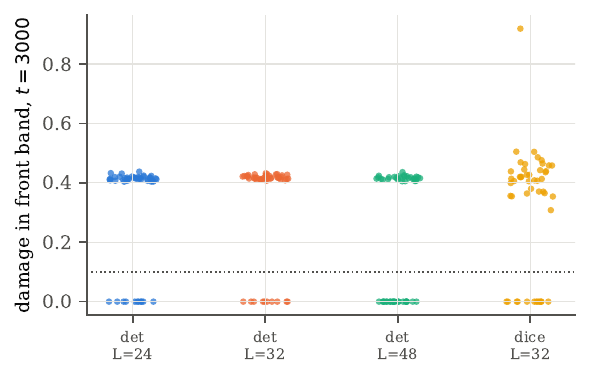}
\caption{E11: damage in the front band (cells within distance 2 of the front), averaged over steps 2801--3000 after a single-unit perturbation (three per front); dotted line: amplification threshold, 0.1.}
\label{fig:front}
\end{figure}

\begin{figure}[t]
\centering
\includegraphics[width=\textwidth]{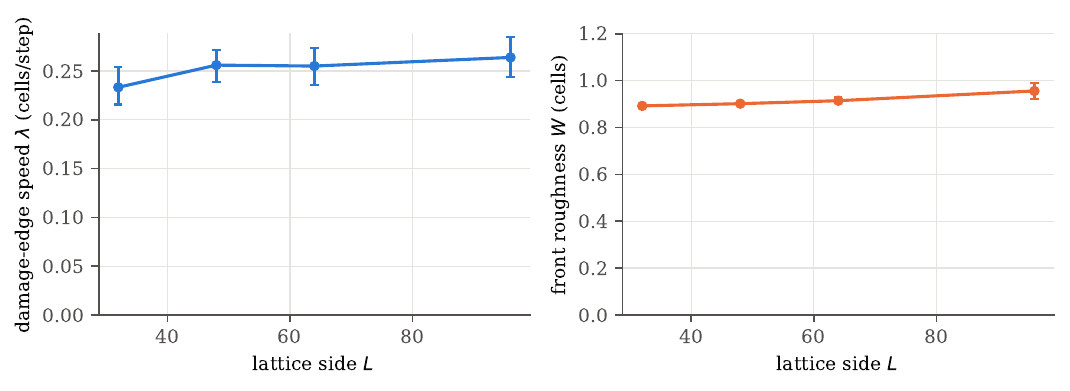}
\caption{E12: speed of the damage edges along the front (left) and front roughness (right) against $L$.}
\label{fig:frontscaling}
\end{figure}

\begin{figure}[t]
\centering
\includegraphics[width=\textwidth]{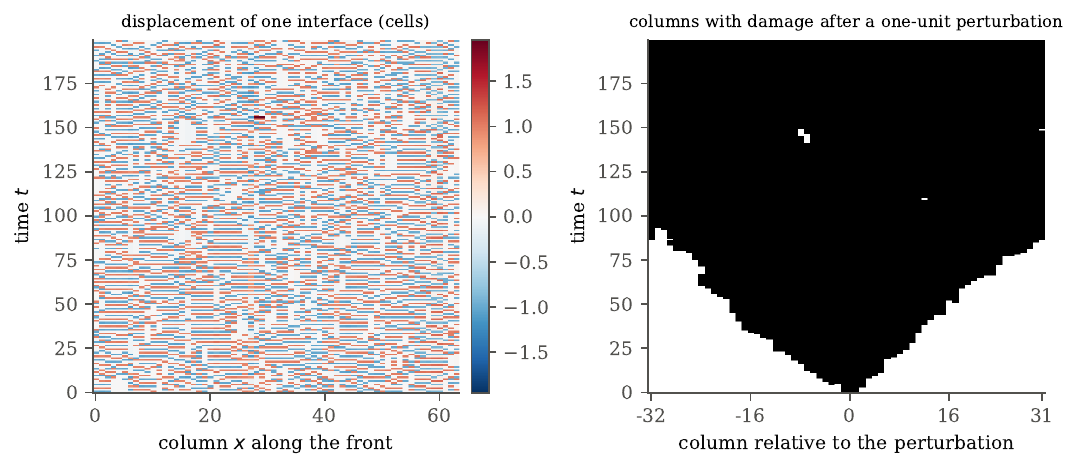}
\caption{Space-time diagrams of one stalemated front ($\mathbb{T}_{64}$, $P=2$, $g=1$, $\theta=1$, deterministic; 200 steps). Left: displacement of one interface in each column, relative to its time mean. Right: columns in which the state differs after a one-unit change of a single front cell at $t=0$; the damage region widens at roughly constant speed until it covers the front.}
\label{fig:frontst}
\end{figure}

\paragraph{Fronts (E11, E12, E14).} With $P=2$, $g=1$ and $\theta=1$, 24--38\% of random starts end in a stalemate (Table~\ref{tab:fronts}), comparable to the probability, about 1/4, with which zero-temperature Ising dynamics on the square lattice freezes into two stripes \cite{spirin2001}; unlike those stripes, RiskCA fronts remain active. Stalemates were generated deliberately (16 per setting at $L=24$, 32, 48 with deterministic battles and at $L=32$ with dice; 8--12 per size at $L=32$ to 96 in E12) and characterised in Tables~\ref{tab:fronts} and \ref{tab:e12}. None of the 92 deterministic fronts repeated a configuration within $2\times10^5$ steps ($10^5$ at $L=96$). Combat continues on a fraction of cells that falls with $L$ in proportion to the fraction on the front, as expected for activity confined to a line. With deterministic battles, 15--30\% of cells change owner at least once within $10^4$ steps but only 3--5\% end with a different owner (41\% and 10\% with dice), so the fronts fluctuate without advancing. Their roughness $W$, the standard deviation over columns $x$ of the stripe thickness $h(x,t)$, averaged over 2000 steps, stays below one cell at all sizes ($W\propto L^{0.06}$). A one-unit change at a front cell is amplified, i.e.\ after 3000 steps it affects more than 10\% of the front band (the cells within distance 2 of the front), in 56--79\% of trials, least often at $L=48$; when amplified, the damage saturates at about 42\% of the band, independently of $L$ and of the battle operator (Figure~\ref{fig:front}). The edges of the damage region move along the front at $\lambda=0.23$ cells per step at $L=32$ and 0.26 at $L=48$--96 (Figure~\ref{fig:frontscaling}), well below the bound of 5 (Lemma~\ref{lem:radius}(b)); the speed rises slightly from $L=32$ to 48 and is constant beyond, as E19 confirms at $L=128$. Pooled over sizes, the left and right speeds are 0.24 and 0.26 (paired Wilcoxon $p=0.05$, Holm $p=0.29$, 107 perturbations), so a small asymmetry cannot be excluded. The one-dimensional space-time pattern $h(x,t)$ of the front is incompressible (compression ratio 0.98--0.99). Figure~\ref{fig:frontst} shows one front: each column moves back and forth by about one cell in an irregular spatial pattern, and a single-unit perturbation spreads along the front as a cone.

\begin{table}[t]
\centering\small
\caption{E14: fronts from six initial conditions ($\mathbb{T}_{48}$, $P=2$, $g=1$, $\theta=1$, deterministic, 12 runs each). Outcome after $2\times10^4$ steps; recurrence within $10^5$ further steps; $\lambda$: damage-edge speed in columns per step; persistent: share of single-unit perturbations not healed after $2L$ steps; excursion: largest distance from the interface at the start of a $10^4$-step window of any cell that changed owner during it.}
\label{tab:geometry}
\begin{adjustbox}{max width=\textwidth}
\begin{tabular}{l c c c c c c c c}\toprule
initial condition & front / cons. / frozen & recurrence & $\lambda$ (columns/step) & persistent & roughness $W$ & compression & activity per front cell & excursion\\ \midrule
random & 3 / 8 / 1 & 0/3 & 0.28 [0.24, 0.31] & 0.50 & 0.90 & 0.99--1.00 & 0.40 [0.39, 0.40] & 2\\
flat-random & 12 / 0 / 0 & 0/12 & 0.26 [0.24, 0.29] & 0.69 & 0.90 & 0.99--1.00 & 0.39 [0.38, 0.40] & 2\\
flat-full & 0 / 0 / 12 & -- & -- & -- & -- & -- & -- & --\\
flat-one & 0 / 0 / 12 & -- & -- & -- & -- & -- & -- & --\\
step & 12 / 0 / 0 & 0/12 & 0.26 [0.24, 0.28] & 0.69 & 0.90 & 0.99--1.00 & 0.39 [0.39, 0.40] & 2\\
diagonal & 12 / 0 / 0 & 0/12 & 0.52 [0.50, 0.55] & 0.77 & 1.09 & 0.99--1.00 & 0.38 [0.38, 0.39] & 2--3\\
\bottomrule\end{tabular}

\end{adjustbox}
\end{table}

E14 tests whether these properties depend on how the front was formed (Table~\ref{tab:geometry}). Prepared straight fronts with random troops, fronts with a step of $L/4$ rows and the stalemates from random starts give the same speed (0.26--0.28 columns per step), roughness (0.90), compression ratio (0.99--1.00) and activity per front cell (0.39--0.40); none recurred within $10^5$ steps. Staircase fronts along the diagonal are also active and non-recurrent, but their damage edges move 0.52 columns per step and their roughness is 1.09, so the speed depends on the orientation of the front relative to the lattice. Straight fronts with uniform troops never fight: a saturated configuration is a fixed point for $\theta\ge1$ (Proposition~\ref{prop:fixed}(b)), and a configuration with one unit per cell grows uniformly, never meets the attack condition, and reaches that fixed point after $\Nmax-1$ steps. Activity requires heterogeneous troops along the front. In every active front, every cell that changed owner within $10^4$ steps lay within 2 cells of the interface at the start of the window (3 for diagonal fronts). For the prepared fronts the two interfaces start $L/2=24$ rows apart, so they stay more than 10 cells apart, as required; for the three stalemates from random starts the separation was not measured. On these windows the conditions stated after Proposition~\ref{prop:local} hold for the prepared straight fronts: each is a one-dimensional CA on $\mathbb{Z}_L$ whose cell is a column segment about ten cells high, and the measured $\lambda$ are finite-window estimates of its left and right Lyapunov exponents in the sense of \citet{shereshevsky1992}; asymptotic exponents would require the reduction to hold indefinitely. For the step and diagonal fronts the same argument applies with the horizontal strip replaced by a band around the interface.

E19 extends the prepared straight fronts to larger lattices and longer times: 8 lattices at $L=64$ followed for $5\times10^5$ steps and 8 at $L=128$ for $2\times10^5$ steps, after a transient of $2\times10^4$ (Table~\ref{tab:longfronts}). All 16 lattices kept both fronts. In every window of $10^4$ steps, every cell that changed owner lay within 2 cells of the front at the start of the long run, the mean position of the fronts moved by less than one row, no configuration repeated, and the activity stayed at 0.39 owner changes per front cell and step. The damage edges moved at 0.26 columns per step at the start and at 0.25--0.29 at the end, at both sizes. Over runs 50 and 20 times longer than the $10^4$-step window used for the excursion in E14, the fronts remain confined and active, and their speed changes neither with $L$ nor with time.

\begin{table}[t]
\centering\small
\caption{E19: prepared straight fronts on larger lattices and over longer times ($P=2$, $g=1$, $\theta=1$, deterministic, 8 lattices per row). Excursion: largest distance, from the front at the start, of a cell that changed owner in any window of $10^4$ steps; drift: largest change of the mean stripe thickness (rows); $\lambda$: damage-edge speed (columns per step) at the start and the end of the run (95\% CI).}
\label{tab:longfronts}
\begin{adjustbox}{max width=\textwidth}
\begin{tabular}{r c c c c c c c}\toprule
$L$ & fronts & steps & excursion (max) & drift (max) & activity per front cell & $\lambda$, start & $\lambda$, end\\ \midrule
64 & 8/8 & $5\times10^{5}$ & 2 & 0.73 & 0.39 [0.39, 0.39] & 0.26 [0.24, 0.28] & 0.29 [0.27, 0.31]\\
128 & 8/8 & $2\times10^{5}$ & 2 & 0.34 & 0.39 [0.39, 0.39] & 0.26 [0.24, 0.28] & 0.25 [0.24, 0.27]\\
\bottomrule\end{tabular}

\end{adjustbox}
\end{table}

\begin{table}[t]
\centering\small
\caption{E11: stalemated fronts ($P=2$, $g=1$, $\theta=1$; 16 per row). Columns give the share of scanned runs ending in stalemate, recurrence, the share of cells on the front, the shares of cells that changed owner at least once and net within $10^4$ steps, and the share of single-unit perturbations amplified after 3000 steps (damage above 10\% of the front band).}
\label{tab:fronts}
\begin{adjustbox}{max width=\textwidth}
\begin{tabular}{l r c c c c c c}\toprule
battle, $L$ & fronts & P(stalemate) & repeated within $2\times10^5$ & front cells & changed owner (10$^4$ steps) & net change & perturbations amplified\\ \midrule
det, 24 & 16 & 0.28 [0.19, 0.40] & 0/16 & 18.4\% & 30.2\% & 5.4\% & 0.77 [0.63, 0.87]\\
det, 32 & 16 & 0.33 [0.23, 0.45] & 0/16 & 13.5\% & 22.9\% & 4.2\% & 0.79 [0.66, 0.88]\\
det, 48 & 16 & 0.38 [0.25, 0.52] & 0/16 & 9.1\% & 15.2\% & 2.8\% & 0.56 [0.42, 0.69]\\
dice, 32 & 16 & 0.24 [0.16, 0.34] & n/a & 14.6\% & 40.6\% & 9.5\% & 0.75 [0.61, 0.85]\\
\bottomrule\end{tabular}

\end{adjustbox}
\end{table}

\begin{table}[t]
\centering\small
\caption{E12: fronts at four sizes (deterministic, four single-unit perturbations per front). $\lambda$: mean speed of the damage edges along the front for spreading perturbations; perturbations that died: healed completely within $2L$ steps; $W$: roughness; last column: compression ratio of the front profile $h(x,t)$.}
\label{tab:e12}
\begin{adjustbox}{max width=\textwidth}
\begin{tabular}{r r c c c c c}\toprule
$L$ & fronts & recurrence & $\lambda$ (cells/step) & perturbations that died & roughness $W$ & compression of $h(x,t)$\\ \midrule
32 & 12 & 0/12 & 0.23 [0.22, 0.25] & 16/48 & 0.89 [0.89, 0.89] & 0.99 [0.99, 0.99]\\
48 & 12 & 0/12 & 0.26 [0.24, 0.27] & 18/48 & 0.90 [0.90, 0.90] & 0.99 [0.99, 1.00]\\
64 & 12 & 0/12 & 0.26 [0.24, 0.27] & 19/48 & 0.91 [0.91, 0.93] & 0.99 [0.99, 0.99]\\
96 & 8 & 0/8 & 0.26 [0.24, 0.28] & 16/32 & 0.96 [0.92, 0.99] & 0.98 [0.97, 0.99]\\
\bottomrule\end{tabular}

\end{adjustbox}
\end{table}

Within the observed windows this one-dimensional CA shows, for all initial geometries with heterogeneous troops, the signatures that operational criteria associate with class III \cite{wolfram1984}: no recurrence, positive Lyapunov speeds that do not change beyond $L=48$, an incompressible space-time pattern, and damage saturation independent of length. No localised propagating structures were observed. We do not claim that the fronts are of class III: non-recurrence only bounds the period of a finite deterministic system from below, a third to a half of single-unit perturbations heal completely within $2L$ steps (E12), and we have no argument, such as a Lyapunov function or a curvature bound, that a front stays in a strip, so the one-dimensional reduction is a finite-window description without a stability guarantee.

\FloatBarrier
\section{Discussion}
\label{sec:discussion}

\paragraph{Mechanisms.} Three mechanisms, each tied to a formal result, account for most of the observations. Combat dissipates troops (Proposition~\ref{prop:balance}), which freezes the rule in finite time without growth (Proposition~\ref{prop:g0}). Saturated cells cannot be attacked when $\theta\ge1$ (condition~\eqref{eq:attack}; Proposition~\ref{prop:fixed}(b)), which freezes fragmented partitions and keeps minority domains alive. Activity is confined to the neighbourhood of interfaces (Proposition~\ref{prop:local}), which makes the fronts the only place where irregular dynamics persists.

\paragraph{Two fields.} The central difference between RiskCA and majority-type lattice models is the second component of the state, $x_i=(o_i,n_i)$. In voter, majority-vote and zero-temperature Ising dynamics \cite{clifford1973,deoliveira1992,spirin2001} the state is the ownership field alone, and whether an interface moves depends only on the local configuration of owners. In RiskCA ownership determines the territorial configuration, but the troop field determines whether a frontier is active: an attack requires $n_i\ge2$ and $n_i-n_j\ge\theta$. Growth raises the troop field, combat depletes it and fortification moves it towards the front, and several of our results follow from this coupling. Saturation freezes partitions at high $\theta$ and leaves minority domains alive in lattices that one faction otherwise dominates (Section~\ref{sec:res:transition}); conversely, interfaces with heterogeneous troops stay active throughout the observed windows, whereas straight Ising stripes are frozen (E11, E14; about 0.4 owner changes per front cell and step in the prepared fronts of E14), and straight fronts with uniform troops never fight. Perturbations of the troop field alone can spread and change ownership (E4b), and most damage changes troop counts without changing territory (E4). A description in terms of ownership alone, as in the models above, would miss the variable that controls activity.

\paragraph{Classes.} Without growth every orbit ends on a fixed point (Proposition~\ref{prop:g0}), which is class-I or class-II behaviour. With growth, runs end on homogeneous attractors at low attack margins and on heterogeneous ones at high margins, saturated fixed points or combinations of fixed local oscillators, after a transient that is sensitive to initial conditions. Since I and II exhaust the attractors by definition, the substantive finding is that long cycles are combinations of fixed local oscillators rather than moving structures (94 of 104 analysed cycles have no cell with the global period), and that oscillators farther apart than the radius showed no interaction in any trial, while contiguous ones interacted in 9 of the 104 cycles (E18), and that structures resembling class IV appear only with index-ordered priorities. This differs from the two-dimensional rules studied by \citet{packard1985}, where all four classes occur on the whole lattice: in RiskCA with growth the irregular dynamics we observed is confined to one-dimensional interfaces, because cells far from every interface change only by growth (Proposition~\ref{prop:local}(a)) and therefore saturate. The labels are insensitive to the thresholds of the classification (Table~\ref{tab:sensitivity}), but they classify attractors reached within $T_{\max}$; 3\% of runs remain unresolved.

\paragraph{Fronts.} The stalemated fronts are the one place where sustained irregular dynamics was observed. Because the dynamics near an interface is locally determined (Proposition~\ref{prop:local}) and the fronts stayed within two or three cells of their initial position (E14), for up to $5\times10^5$ steps and on lattices up to $128\times128$ (E19), each prepared straight front is, over the observed windows, a one-dimensional CA of finite radius with the operational signatures of class III. For the reasons given in Section~\ref{sec:res:robust}, this is a finite-window description, not a class-III classification. The fronts are also anisotropic: damage moves twice as fast, in columns per step, along diagonal fronts as along straight ones (0.52 against 0.26), so the front dynamics depends on the orientation of the interface relative to the lattice.

\paragraph{Battle operator.} Over the $(g,\theta)$ plane, the dice reproduce the regime structure of the deterministic rule (E16, Table~\ref{tab:dicescan}): consolidation at low $\theta$, with probabilities differing by at most 0.19 between the operators on the same initial states, and saturated multi-owner partitions for $\theta\ge4$ (94--100\% of dice runs). The majority outcome, consolidation or a multi-owner state, agrees in all 16 parameter cells, and for both operators it changes between $\theta=2$ and 3. Table~\ref{tab:detdice} collects the other observables measured with both operators at $g=1$, $\theta=1$. The exact dice and the deterministic reduction give the same scaling of consolidation time, the same contagion and coarsening exponents, and the same frequency of active fronts on the torus. Dice slow the spreading of damage and make fronts more mobile (more owner changes, larger net displacement). The lower stalemate probability with dice (0.24 against 0.33, Fisher $p=0.26$) and the larger number of active fronts on open boundaries (13 against 6 of 24, Fisher $p=0.08$) are not significant at the 5\% level. $\mathcal{R}_\rho$ is therefore adequate for the regime structure, which is what the classification uses, but not a substitute for the dice in quantitative statements about damage or fronts, and the comparison says nothing about individual trajectories, which differ between the operators. The damage comparison was made at $g=1$, $\theta=1$ only; whether the discrepancy persists at other parameters was not tested.

\paragraph{Spatial structure.} With identical rules, the time to consolidation grows as $N^{1.0}$--$N^{1.1}$ on the lattice but is nearly independent of $N$ both in the mean-field control and on random 4-regular graphs with local interactions (E5, E5b), and conquests are clustered on the lattice ($\kappa\approx5$) but hardly on the random graph ($\kappa\approx1.5$). This points to the low-dimensional geometry of the lattice as the cause, rather than the locality of interaction, with the caveat that the random graph is locally tree-like (Section~\ref{sec:limits}). The coarsening, the $L^2$ transients and the stripe states are shared with majority-type and zero-temperature Ising dynamics \cite{bray1994,spirin2001}.

\paragraph{Methodology.} Index-ordered tie-breaking produces translating structures that mimic class IV (E2, E3), and a window of 8000 steps would have left long transients unresolved (E10b). Classifications of game-derived CA therefore need a symmetric tie-breaking control and observation windows long enough to exclude long transients.

Hypothesis~\ref{hyp:edge} is supported only in a weak sense: transients are longest at intermediate $\theta$ (E2, E13), but the attractors are fixed points on both sides of the dominance--fragmentation crossover, and sustained irregular dynamics appears only on fronts, not on whole lattices at intermediate $\theta$. In the sense of \citet{langton1990}, in which complex dynamics persists near the boundary between ordered and chaotic regimes, the hypothesis is not supported.

\section{Limitations}
\label{sec:limits}

\paragraph{Model.} RiskCA omits much of the game: cards, discretionary retreat, deliberate reinforcement and strategic play. All cells follow one fixed heuristic, fortification is short-ranged, and the cap $\Nmax$ has no counterpart in the game. Only the battle kernel is validated exactly; this does not validate the attack, growth and fortification protocol, and the board statistics are compared only with qualitative features of the game, for lack of data on human play. No interpretation of $\theta$, $g$ or the outcomes in terms of real conflicts has been validated: RiskCA is a theoretical study of a local rule, not a predictive model of games or conflicts. Several results, notably $L^2$ coarsening, curvature-driven interface motion and stripe states, are generic to majority-type lattice models rather than specific to \emph{Risk}.

\paragraph{Classification.} Wolfram's classes are asymptotic and undecidable in general \cite{culik1988}, and our criteria are operational. Attractor classes are assigned only to runs that reach an attractor within $T_{\max}$; unresolved runs are labelled by finite-time damage, which measures the transient and not the attractor. The thresholds were fixed after a pilot sweep; varying them does not change the majority labels (Table~\ref{tab:sensitivity}), but other indicators could.

\paragraph{Statistics.} The scaling exponents rest on five sizes with 12 runs each and a large run-to-run scatter. They are a consistency check with $L^2$ for runs that reach an attractor; the censored runs, many of them stalemates, form a second population that a single censored model describes poorly, and a competing-risks analysis would need more runs. The dominance--fragmentation step is resolved with a continuous threshold only up to $L=96$, with 24 runs per point and some runs censored, which is not enough to decide whether it becomes a phase transition. The fronts are described over finite windows (at most $5\times10^5$ steps) with 8--16 fronts per setting and lattices of at most $128\times128$, and their confinement is observed, not proved.

\paragraph{Controls.} The mean-field control changes both geometry and mixing; the random-regular-graph control separates the two but has its own structure (short cycles are rare, so it is locally tree-like). Robustness was tested for the battle operator, the fortification fraction, tie-breaking, the updating scheme, boundaries and neighbourhood, each separately and mostly at $L=32$; the order of the sub-rules and combinations of variations were not tested.

\section{Conclusion}
\label{sec:conclusion}

The attack and fortification rules of \emph{Risk}, formulated as a cellular automaton, define a rule of radius at most 6 (at most 5 for deterministic battles with $\rho>1$ and $\theta\ge1$), both bounds being attained, that, for the default parameters, is irreversible in its deterministic form, and that freezes in finite time without growth. With growth, increasing the attack margin moves the system from dominance by one faction to a frozen fragmented partition. The largest faction's share drops between the same two threshold values at all sizes up to $L=128$, more steeply in larger systems, while the probability of complete consolidation is a size-dependent crossover driven by saturated minority domains. Compared with a mean-field control or a random regular graph with the same degree, the lattice lengthens conflicts about twentyfold and makes conquests spatially clustered, which points to its low-dimensional geometry. Whole lattices end on homogeneous or saturated attractors or on long cycles that combine spatially fixed local oscillators; with quenched priorities, no moving structures were found. Straight fronts spanning the lattice stayed within two cells of their initial position for up to $5\times10^5$ steps and on lattices up to $128\times128$, and over these windows they are one-dimensional CA with the operational signatures of class III; whether they remain confined indefinitely is open. Structures suggestive of class IV arise only with index-ordered tie-breaking.

Three questions remain open. First, whether the step in the largest faction's share becomes a discontinuity for $L\to\infty$: with a continuous threshold it is about half a unit of $\theta$ wide up to $L=96$, and deciding requires larger systems and more runs. Second, whether the fronts stay confined and active indefinitely, which would need an argument such as a Lyapunov function or a curvature bound, and whether they are more or less irregular with dice, under which they are more mobile (Table~\ref{tab:detdice}). Third, how the results change with heterogeneous terrain, strategically chosen attack rules and irregular boards.

\paragraph{Data and code availability.} Code, data, the scripts that generate all figures and tables, and the Lean formalisation are available at \url{https://github.com/MaTeArNog/RiskAutomata}.

\paragraph{Funding.} This research received no external funding.

\paragraph{Declaration of interests.} None.

\paragraph{Use of AI tools.} Claude (Anthropic) was used to write the code of the experiments and the Lean formalisation, and to assist in drafting and editing the text. The author directed the work, checked the results and takes full responsibility for the content.

\newcommand{\doi}[1]{\href{https://doi.org/#1}{doi:#1}}

\FloatBarrier
\appendix
\section{Notation}
\label{app:notation}

\begin{table}[h]
\centering\small
\caption{Main symbols and labels.}
\label{tab:notation}
\begin{tabular}{ll}
\toprule
Symbol & Meaning \\
\midrule
$N$, $L$, $\mathbb{T}_L$ & number of cells; side of the square lattice; $L\times L$ torus \\
$P$, $\Nmax$ & number of factions; troop cap per cell \\
$o_i$, $n_i$ & owner and troops of cell $i$ \\
$g$, $\theta$ & growth per step; attack margin \\
$\theta_c$ & effective attack threshold of the stochastic variant E17 \\
$\mathcal{R}_\rho$, $\rho$ & deterministic battle operator; exchange ratio (1.1718) \\
$\mathcal{R}_{\mathrm{dice}}$ & battle operator with the exact dice distribution \\
$\pi$ & tie-breaking priority (quenched random, or index order) \\
$\eta_i$, $\delta_i$ & number of hostile neighbours; distance to the nearest frontier cell \\
$\tau$, $p$, $T_{\max}$ & transient; period; observation window \\
I, II & attractor with one owner; with several owners \\
U$_\mathrm{s}$, U$_\mathrm{l}$ & no attractor within $T_{\max}$, with spreading or localised owner damage \\
$H_o$, $H_f$ & fraction of cells with different owner; with different owner or troops \\
$C$, $\sigma_S$ & compression ratio; standard deviation of the event entropy \\
$\kappa$ & contagion ratio of conquests \\
$\lambda$, $W$ & speed of the damage edges along a front; front roughness \\
\bottomrule
\end{tabular}
\end{table}

\section{CA pseudocode}
\label{app:code}

Algorithm~\ref{alg:step} is one synchronous step of RiskCA as implemented in \texttt{riskca/model.py}, and Algorithm~\ref{alg:kernel} computes the battle kernel. All operations act on all cells at once, and $\arg\min$ and $\arg\max$ over pairs compare lexicographically.

\begin{algorithm}[h]
\caption{One synchronous RiskCA step $x\mapsto\Phi(x)$}
\label{alg:step}
\small
\begin{algorithmic}[1]
\Require owners $o\in\{1..P\}^N$, troops $n\in\{1..\Nmax\}^N$, neighbour table $\mathcal{N}$, priority $\pi$, parameters $(g,\theta,\rho,\Nmax)$
\State \textbf{(G)} $n_i\gets\min(\Nmax,n_i+g)$ for all $i$ \Comment{Risk-like reinforcement on the 42-territory board}
\State \textbf{(D)} \textbf{for all} $i$: $H_i\gets\{j\in\mathcal{N}(i):o_j\neq o_i\}$
\State \hspace{1.2em}\textbf{if} $H_i\neq\emptyset$: $t_i\gets\arg\min_{j\in H_i}(n_j,\pi(j))$;\ $w_i\gets[n_i\ge2\wedge n_i-n_{t_i}\ge\theta]$ \textbf{else} $w_i\gets$ false
\State \textbf{(R)} \textbf{for all} $j$ with $A_j=\{i:w_i\wedge t_i=j\}\neq\emptyset$: $s(j)\gets\arg\max_{i\in A_j}(n_i-1,-\pi(i))$ \Comment{one attacker per target; no attack if $A_j=\emptyset$}
\State \hspace{1.2em}$S\gets\{s(j)\}$; \textbf{for} $i\in S$: $a_i\gets n_i-1$, $n_i\gets1$ \Comment{commit; origin keeps one unit}
\State $\hat o\gets o$ \Comment{owners at the start of (R)}
\For{$i\in S$ \textbf{(in parallel)}, $j=t_i$, $d=n_j$}
  \If{deterministic}
    \If{$\rho a_i>d$} $o_j\gets \hat o_i$; $n_j\gets\max(1,a_i-\lceil d/\rho\rceil)$
    \Else\ $n_j\gets\max(1,d-\lfloor\rho a_i\rfloor)$
    \EndIf
  \Else\ draw $k\sim K(a_i,d;\cdot)$ (Alg.~\ref{alg:kernel}); \textbf{if} attacker survives with $r\ge1$: $o_j\gets \hat o_i$, $n_j\gets r$ \textbf{else} $n_j\gets$ defender survivors
  \EndIf
\EndFor
\State \textbf{(F)} recompute $H_i$ and $\eta_i=|H_i|$ on the new owners
\State \hspace{1.2em}\textbf{for all} $i$ with $\eta_i=0$: $k\gets\arg\max_{k\in\mathcal{N}(i),o_k=o_i}(\eta_k,-\pi(k))$; $m_i\gets\lfloor(n_i-n_k)/2\rfloor$
\State \hspace{1.2em}\textbf{if} $\eta_k>0\wedge m_i\ge1$: $n_i\gets n_i-m_i$; inflow$_k\mathrel{+}=m_i$
\State \hspace{1.2em}$n\gets\min(\Nmax,n+\text{inflow})$ \Comment{overflow is lost}
\State \Return $(o,n)$
\end{algorithmic}
\end{algorithm}

\begin{algorithm}[h]
\caption{Exact battle kernel $K(a,d;\cdot)$ (Tan--Osborne)}
\label{alg:kernel}
\small
\begin{algorithmic}[1]
\State \textbf{for} $\alpha\in\{1,2,3\},\delta\in\{1,2\}$: enumerate all $6^{\alpha+\delta}$ rolls; sort both sides in decreasing order; compare pairwise (ties to defender); tabulate $q_{\alpha\delta}(\ell_A,\ell_D)$
\State $K(a,0)\gets\delta_{\text{att wins},a}$;\ $K(0,d)\gets\delta_{\text{def holds},d}$
\For{$a=1..A$, $d=1..D$ (increasing)}
  \State $K(a,d)\gets\sum_{(\ell_A,\ell_D)}q_{\min(3,a)\min(2,d)}(\ell_A,\ell_D)\,K(a-\ell_A,d-\ell_D)$
\EndFor
\end{algorithmic}
\end{algorithm}

\paragraph{Code layout.} \texttt{riskca/dice.py} (kernel), \texttt{riskca/boards.py} (torus, 42-territory board, mean-field neighbours, random regular graphs), \texttt{riskca/model.py} (rule, schedules, recording), \texttt{riskca/metrics.py} (damage, compression, bootstrap/Wilson CIs), \texttt{experiments/e1--e16*.py} (one script per experiment), \texttt{experiments/configs.py} (machine-readable list of configurations), \texttt{tests/} (unit tests), \texttt{experiments/analyse.py} (all statistics, tables and figures), \texttt{experiments/run\_all.sh} (end-to-end reproduction).

\section{A witness for radius 6}
\label{app:witness}
Table~\ref{tab:witness} gives the configuration $X$ on a $6\times5$ patch with open boundaries ($P=2$, $g=1$, $\theta=0$, $\Nmax=16$, deterministic battles, priority $\pi$ equal to the row-major index, so that cell $(r,c)$ has priority $5(r-1)+c-1$ and lower values win ties). The configuration $Y$ differs from $X$ only at $u=(6,5)$, whose owner is changed from $B$ to $A$. After one step the two images differ at $v=(3,2)$, at graph distance $3+3=6$ from $u$. The difference propagates as follows; all quantities are after growth.
\begin{enumerate}[leftmargin=*,itemsep=0pt]
\item Distance 1. In $X$, cell $(6,4)$ ($A$, 16 units) attacks its weakest hostile neighbour $(6,5)$ (2 units). In $Y$, $(6,5)$ is friendly and $(6,4)$ attacks $(6,3)$ ($B$, 5 units) instead.
\item Distance 2. In $X$, the only attack on $(6,3)$ comes from $(5,3)$ ($A$, 5 units, $a=4$) and is accepted. In $Y$, $(6,4)$ commits $a=15>4$, so its attack on $(6,3)$ is accepted and that of $(5,3)$ is not.
\item Distance 3. In both configurations $(6,3)$ attacks $(5,3)$ with $a=4$. In $X$, $(5,3)$ has committed its own attack and defends with one unit, so it is conquered by $B$. In $Y$, it defends with 5 units, and $\rho\cdot4=4.69<5$, so it holds.
\item Distance 4. After (R), cell $(4,3)$ ($A$) has two hostile neighbours, $(4,4)$ and $(5,3)$, in $X$, and one, $(4,4)$, in $Y$.
\item Distance 5. Cell $(3,3)$, conquered by $(4,3)$ with 14 units, is interior after (R). Its friendly neighbours of highest threat are $(4,3)$ alone in $X$ ($\eta=2$), and $(4,3)$ and $(3,2)$ in $Y$ ($\eta=1$ each, since $(3,1)$ has been conquered by $B$). In $Y$ the tie goes to $(3,2)$, which has the lower priority value.
\item Distance 6. Cell $(3,3)$ sends 6 units to $(4,3)$ in $X$ and to $(3,2)$ in $Y$, so $(3,2)$ ends with 2 units in $X$ and 8 in $Y$.
\end{enumerate}
The example is reproduced by a unit test of the model code. Along the six steps, $\pi$ is needed only in step 5, where it breaks the tie. The row-major order used here is the index order that produces artefacts in Section~\ref{sec:res:artefact}; it is used only to make the example explicit. The bound of Lemma~\ref{lem:radius} holds for every $\pi$, and the random search from which this example was reduced (two cases among $6\times10^4$ single-cell perturbations on $\mathbb{T}_{16}$) used random quenched priorities, so distance 6 is reached without index order. Other ties in the example also depend on $\pi$, so the reduced configuration itself is specific to the order given; with $\theta\ge1$ and deterministic battles step 3 is impossible, because every accepted attack succeeds (Lemma~\ref{lem:radius}(b)).

\begin{table}[h]
\centering\small
\caption{Radius witness: owner and troops of each cell of $X$ (rows 1--6 from the top, columns 1--5 from the left). Bold: the perturbed cell $u$; underlined: the cell $v$ at distance 6 whose image changes.}
\label{tab:witness}
\begin{tabular}{c ccccc}
\toprule
 & 1 & 2 & 3 & 4 & 5\\
\midrule
1 & A1 & A1 & A1 & A1 & A1\\
2 & A1 & A1 & A1 & A16 & A1\\
3 & A1 & \underline{A1} & B1 & B1 & A1\\
4 & B1 & A1 & A16 & B16 & A1\\
5 & A1 & A1 & A4 & A1 & A1\\
6 & A1 & B1 & B4 & A16 & \textbf{B1}$\to$\textbf{A1}\\
\bottomrule
\end{tabular}
\end{table}

\section{Machine-checked results}
\label{app:lean}
The Lean~4 formalisation (directory \texttt{lean/} of the repository) implements the rule independently of the Python code, with the battle operator as a parameter. Table~\ref{tab:lean} lists the main declarations; each is checked by the Lean kernel without \texttt{sorry} or \texttt{native\_decide}. ``Every realisation'' means every battle outcome function in which a winning attacker keeps between 1 and $a$ units and a holding defender between 1 and $d$; every outcome of the dice, and the dice operator of the code for every $\Nmax$ and every choice of random numbers, are of this kind (\texttt{diceReal\_ok}, \texttt{diceUN\_ok}). Agreement with the Python implementation is checked on concrete configurations: the Lean model reproduces the outputs of the code on the example of Figure~\ref{fig:schematic}, the three radius witnesses and the 19 dice battles of the dice witness. Boards are finite graphs; Propositions~\ref{prop:irrev} and \ref{prop:fixed}(c) assume connected boards, as stated. The experimental results of Sections~\ref{sec:results} and \ref{sec:discussion} are not formalised.

\begin{table}[h]
\centering\small
\caption{Lean declarations for the formal results.}
\label{tab:lean}
\begin{tabular}{@{}>{\raggedright\arraybackslash}p{4.0cm}>{\raggedright\arraybackslash}p{6.6cm}>{\raggedright\arraybackslash}p{3.5cm}@{}}
\toprule
Result & Lean declaration & Battle operator \\
\midrule
Lemma~\ref{lem:radius}(a), (c) & \texttt{\footnotesize lemma1a}, \texttt{\footnotesize lemma1c} & every realisation \\
Lemma~\ref{lem:radius}(b) & \texttt{\footnotesize lemma1b} & deterministic \\
Lemma~\ref{lem:radius}, bounds attained & \texttt{\footnotesize w6\_reach}, \texttt{\footnotesize w5\_reach}, \texttt{\footnotesize w6d\_reach} & deterministic; dice \\
Lemma~\ref{lem:radius}, spread per step & \texttt{\footnotesize light\_cone} & every realisation \\
Proposition~\ref{prop:finite} & \texttt{\footnotesize prop1} & every fixed realisation \\
Proposition~\ref{prop:irrev} & \texttt{\footnotesize prop2}, \texttt{\footnotesize prop2\_garden}, \texttt{\footnotesize prop2\_identity}, \texttt{\footnotesize prop2\_cases} & deterministic \\
Proposition~\ref{prop:balance} & \texttt{\footnotesize prop3}, \texttt{\footnotesize prop3\_terms} & every realisation \\
Proposition~\ref{prop:fixed}(a), (b) & \texttt{\footnotesize prop4a\_homogeneous}, \texttt{\footnotesize prop4a\_time}, \texttt{\footnotesize prop4b\_saturated\_fixed}, \texttt{\footnotesize prop4b\_count} & every realisation \\
Proposition~\ref{prop:fixed}(c) & \texttt{\footnotesize prop4c} & deterministic \\
Proposition~\ref{prop:g0} & \texttt{\footnotesize prop5} & every sequence of realisations \\
Proposition~\ref{prop:local} & \texttt{\footnotesize prop6a}, \texttt{\footnotesize prop6b} & every realisation \\
Figure~\ref{fig:schematic}, Appendix~\ref{app:witness} & \texttt{\footnotesize fig1}, \texttt{\footnotesize w6\_story} & deterministic \\
Section~\ref{sec:res:dice}: rolls, $\rho$, kernel, $\mathcal{R}_\rho$ errors & \texttt{\footnotesize roll\_3v2}, \texttt{\footnotesize exchange\_ratio}, \texttt{\footnotesize win\_probs}, \texttt{\footnotesize rho\_errors}, \texttt{\footnotesize rho\_survivor\_error} & dice \\
\bottomrule
\end{tabular}
\end{table}

\section{The 42-territory board}
\label{app:board}
The six regions (bonus; territories) are: N.\ America (5; Alaska, Northwest Territory, Greenland, Alberta, Ontario, Quebec, Western US, Eastern US, Central America); S.\ America (2; Venezuela, Peru, Brazil, Argentina); Europe (5; Iceland, Scandinavia, Ukraine, Great Britain, Northern Europe, Western Europe, Southern Europe); Africa (3; North Africa, Egypt, East Africa, Congo, South Africa, Madagascar); Asia (7; Ural, Siberia, Yakutsk, Kamchatka, Irkutsk, Mongolia, Japan, Afghanistan, China, Middle East, India, Siam); Australia (2; Indonesia, New Guinea, Western Australia, Eastern Australia). The 83 edges are listed in \texttt{riskca/boards.py}. They include the inter-regional links Alaska--Kamchatka, Greenland--Iceland, Central America--Venezuela, Brazil--North Africa, Western Europe--North Africa, Southern Europe--\{North Africa, Egypt, Middle East\}, Ukraine--\{Ural, Afghanistan, Middle East\}, Egypt--Middle East, East Africa--Middle East and Siam--Indonesia. The number of \emph{entry cells} of a region (cells with at least one neighbour outside it) is 1 (Australia), 2 (S.\ America), 3 (N.\ America, Africa), 4 (Europe) and 5 (Asia).

\section{Supplementary tables}
\label{app:supp}
Tables~\ref{tab:sensitivity}--\ref{tab:dicescan} give the details of the robustness and sensitivity analyses summarised in Sections~\ref{sec:results} and \ref{sec:discussion}.

\begin{table}[htbp]
\centering\small
\caption{Sensitivity of the E2 labels to the thresholds (quenched tie-breaking, 640 runs, 40 parameter cells). $T_{\max}$ variants re-censor runs whose attractor was detected only after the shorter window.}
\label{tab:sensitivity}
\begin{adjustbox}{max width=\textwidth}
\begin{tabular}{l c c c}\toprule
variant & runs relabelled & cells with new majority label & unresolved runs\\ \midrule
default ($H_o$ threshold 0.4, $T_{\max}=3\times10^4$) & -- & -- & 20\\
$H_o$ threshold 0.2 & 2/640 & 0/40 & 20\\
$H_o$ threshold 0.3 & 1/640 & 0/40 & 20\\
$H_o$ threshold 0.5 & 3/640 & 0/40 & 20\\
$H_o$ threshold 0.6 & 4/640 & 0/40 & 20\\
$T_{\max}=8000$ & 9/640 & 0/40 & 29\\
$T_{\max}=3000$ & 17/640 & 0/40 & 37\\
$T_{\max}=1000$ & 53/640 & 1/40 & 73\\
\bottomrule\end{tabular}

\end{adjustbox}
\end{table}

\begin{table}[htbp]
\centering\small
\caption{E4b: distribution of damage at $t=600$ over 64 pairs per row ($\mathbb{T}_{64}$, deterministic). Perturbation: change of owner, or one unit added (removed at the cap). Healed: $H_f(600)=0$. Last column: largest absolute Spearman correlation of $H_o(600)$ with the troops or the number of hostile neighbours of the perturbed cell.}
\label{tab:damagehet}
\begin{adjustbox}{max width=\textwidth}
\begin{tabular}{l l c c c c}\toprule
regime & perturbation & $H_o(600)$: 10\%, 50\%, 90\% & $H_f(600)$: 50\% & healed & max $|\rho_s|$\\ \midrule
R2 ($g$=1, $\theta$=1) & owner & 0.000, 0.049, 0.159 & 0.243 & 0.19 & 0.15\\
R2 ($g$=1, $\theta$=1) & one unit & 0.000, 0.035, 0.137 & 0.208 & 0.39 & 0.23\\
R3 ($g$=2, $\theta$=4) & owner & 0.000, 0.136, 0.201 & 0.215 & 0.20 & 0.34\\
R3 ($g$=2, $\theta$=4) & one unit & 0.000, 0.113, 0.160 & 0.172 & 0.25 & 0.29\\
R4 ($g$=1, $\theta$=8) & owner & 0.517, 0.648, 0.686 & 0.668 & 0.06 & 0.12\\
R4 ($g$=1, $\theta$=8) & one unit & 0.000, 0.000, 0.663 & 0.000 & 0.61 & 0.16\\
\bottomrule\end{tabular}

\end{adjustbox}
\end{table}

\begin{table}[htbp]
\centering\small
\caption{PRCC values; $^*$: 95\% CI excludes zero.}
\label{tab:prcc}
\begin{tabular}{lcccccc}\toprule
output & $g$ & $\theta/N_{max}$ & $N_{max}$ & $\rho$ & $P$ & $\mathrm{fortify}$\\ \midrule
log10 transient & +0.34$^*$ & -0.48$^*$ & +0.27$^*$ & +0.02 & +0.22$^*$ & -0.05\\
consolidated & +0.09$^*$ & -0.46$^*$ & -0.01 & -0.01 & -0.04 & -0.10$^*$\\
final interface & +0.09 & +0.81$^*$ & +0.02 & +0.05 & +0.55$^*$ & +0.04\\
mean activity & -0.41$^*$ & -0.87$^*$ & +0.23$^*$ & +0.03 & +0.35$^*$ & -0.02\\
\bottomrule\end{tabular}

\end{table}

\begin{table}[htbp]
\centering\small
\caption{E10: majority label with the Moore neighbourhood ($L=32$); $^\dagger$: differs from the von Neumann result.}
\label{tab:moore}
\begin{tabular}{ccccccccc}\toprule
$g\backslash\theta$ & 0 & 1 & 2 & 3 & 4 & 6 & 8 & 12\\ \midrule
0 & II & II & II & II & II & II & II & II\\
1 & I & I & I & I$^\dagger$ & II & II & II & II\\
2 & I & I & I & I$^\dagger$ & II & II & II & II\\
\bottomrule\end{tabular}

\end{table}

\begin{table}[htbp]
\centering\small
\caption{E15: classification with fortification fraction $f$ ($\mathbb{T}_{32}$, $g\in\{1,2\}$, the eight values of $\theta$ of E2, 16 runs each, same seeds as E2).}
\label{tab:fort}
\begin{adjustbox}{max width=\textwidth}
\begin{tabular}{l c c c c c}\toprule
fraction $f$ & cells with new majority label & P(consolidation) & multi-owner fixed points (saturated) & long cycles & unresolved\\ \midrule
1/2 (E2) & -- & 0.35 & 71 (71) & 53 & 12\\
1/3 & 2/16 & 0.34 & 82 (82) & 52 & 17\\
2/3 & 0/16 & 0.35 & 60 (60) & 58 & 11\\
\bottomrule\end{tabular}

\end{adjustbox}
\end{table}

\begin{table}[htbp]
\centering\small
\caption{Deterministic ($\mathcal{R}_\rho$) and exact dice ($\mathcal{R}_{\mathrm{dice}}$) battles on the same observables ($g=1$, $\theta=1$).}
\label{tab:detdice}
\begin{adjustbox}{max width=\textwidth}
\begin{tabular}{l l c c}
\toprule
Observable & Experiment & $\mathcal{R}_\rho$ & $\mathcal{R}_{\mathrm{dice}}$ \\
\midrule
Consolidation time exponent $z_N$ & E5 & 1.09 [0.95, 1.22] & 1.03 [0.57, 1.44] \\
Median consolidation time, $N=2304$ & E5 & 767 & 934 \\
Contagion $\kappa$ & E5 & 5.2 [4.9, 5.6] & 4.9 [4.7, 5.2] \\
Interface decay exponent (window 20--1000) & E5 & 0.67 & 0.67 \\
Damage-front speed (cells/step) & E4 & 0.53 & 0.37 \\
Owner / full-state damage at $t=600$ & E4 & 0.12 / 0.31 & 0.05 / 0.17 \\
Active fronts, torus ($P=2$; $P=4$) & E9 & 7/24; 3/24 & 7/24; 3/24 \\
Active fronts, open boundaries ($P=2$; $P=4$) & E9 & 6/24; 3/24 & 13/24; 4/24 \\
P(stalemate), $P=2$, $L=32$ & E11 & 0.33 [0.23, 0.45] & 0.24 [0.16, 0.34] \\
Front cells changing owner / net, $10^4$ steps & E11 & 22.9\% / 4.2\% & 40.6\% / 9.5\% \\
Perturbations amplified at fronts & E11 & 0.79 [0.66, 0.88] & 0.75 [0.61, 0.85] \\
\bottomrule
\end{tabular}
\end{adjustbox}
\end{table}

\begin{table}[htbp]
\centering\small
\caption{E16: deterministic and dice battles from the same initial states ($\mathbb{T}_{32}$, 16 runs per row). Active: not consolidated, with owners still changing at $T_{\max}=3\times10^4$ (deterministic: no attractor found; dice: more than $10^{-3}$ of cells change owner per step over the last 500 steps). Last column: share of dice runs ending in a saturated multi-owner configuration.}
\label{tab:dicescan}
\begin{adjustbox}{max width=\textwidth}
\begin{tabular}{c c cc cc c}\toprule
 & & \multicolumn{2}{c}{P(consolidation)} & \multicolumn{2}{c}{P(active at $T_{\max}$)} & dice: saturated\\
$g$ & $\theta$ & det. & dice & det. & dice & multi-owner\\ \midrule
1 & 0 & 0.75 & 0.94 & 0.19 & 0.06 & 0.00\\
1 & 1 & 0.88 & 0.94 & 0.12 & 0.06 & 0.00\\
1 & 2 & 0.62 & 0.75 & 0.19 & 0.12 & 0.12\\
1 & 3 & 0.19 & 0.25 & 0.00 & 0.12 & 0.62\\
1 & 4 & 0.06 & 0.06 & 0.00 & 0.00 & 0.94\\
1 & 6 & 0.00 & 0.00 & 0.00 & 0.00 & 1.00\\
1 & 8 & 0.00 & 0.00 & 0.06 & 0.00 & 1.00\\
1 & 12 & 0.00 & 0.00 & 0.00 & 0.00 & 1.00\\
2 & 0 & 1.00 & 1.00 & 0.00 & 0.00 & 0.00\\
2 & 1 & 0.88 & 0.94 & 0.12 & 0.06 & 0.00\\
2 & 2 & 0.81 & 0.81 & 0.06 & 0.12 & 0.06\\
2 & 3 & 0.44 & 0.38 & 0.00 & 0.00 & 0.62\\
2 & 4 & 0.00 & 0.00 & 0.00 & 0.00 & 1.00\\
2 & 6 & 0.00 & 0.00 & 0.00 & 0.00 & 1.00\\
2 & 8 & 0.00 & 0.00 & 0.00 & 0.00 & 1.00\\
2 & 12 & 0.00 & 0.00 & 0.00 & 0.00 & 1.00\\
\bottomrule\end{tabular}

\end{adjustbox}
\end{table}

\FloatBarrier

\section{Hypothesis tests}
\label{app:tests}
Table~\ref{tab:holm} lists the hypothesis tests reported in the paper with raw and Holm-adjusted $p$-values.

\begin{table}[h]
\centering\small
\caption{Hypothesis tests reported in the paper, with Holm-adjusted $p$-values.}
\label{tab:holm}
\begin{tabular}{l r r}\toprule
test & $p$ & Holm-adjusted $p$\\ \midrule
E2 Mann--Whitney, $H_o$, class I vs.\ II & $<10^{-17}$ & $<10^{-16}$\\
E2 Spearman, $H_o$ vs.\ $C$ & 0.35 & 0.79\\
E5 Mann--Whitney, lattice vs.\ mean-field (dice, $L=32$) & $<10^{-6}$ & $<10^{-5}$\\
E6 Kruskal--Wallis, $T_c$ across $P$ (det) & 0.0059 & 0.047\\
E6 Spearman, $T_c$ vs.\ $P$ (det) & 0.053 & 0.29\\
E6 Kruskal--Wallis, $T_c$ across $P$ (dice) & 0.0018 & 0.016\\
E6 Spearman, $T_c$ vs.\ $P$ (dice) & 0.00063 & 0.0067\\
E7 Mann--Whitney, sync vs.\ turn ($P=4$) & 0.00061 & 0.0067\\
E7 Mann--Whitney, map vs.\ mean-field (turn, $P=4$) & $<10^{-132}$ & $<10^{-131}$\\
E7 Kruskal--Wallis, length across $P$ (turn) & $<10^{-196}$ & $<10^{-195}$\\
E7 Kruskal--Wallis, length across $P$ (sync) & $<10^{-237}$ & $<10^{-235}$\\
E7 Spearman, holding vs.\ region size & 0.00031 & 0.0037\\
E7 Spearman, holding vs.\ entry cells & 0.015 & 0.1\\
E9 Fisher, active fronts, periodic vs.\ open & 0.4 & 0.79\\
E9 Fisher, active fronts, open, det.\ vs.\ dice ($P=2$) & 0.075 & 0.3\\
E11 Fisher, stalemate, det.\ vs.\ dice ($L=32$) & 0.26 & 0.79\\
E12 Wilcoxon, left vs.\ right damage-edge speed & 0.048 & 0.29\\
\bottomrule\end{tabular}

\end{table}

\end{document}